\documentclass[review]{elsarticle}

\usepackage{lineno,hyperref,bm,psfrag,graphicx,amsmath,amssymb,subfigure,multirow,graphpap,hyperref,color}
\modulolinenumbers[5]

\newtheorem{remark}{Remark}[section]

\newcommand{\aver}[1]{\left\{\!\!\left\{#1\right\}\!\!\right\}}

\newcommand{\jump}[1]{\left[\!\left[#1\right]\!\right]}
\newcommand{\bigjump}[1]{\left[\!\!\left[#1\right]\!\!\right]}

\journal{Computers $\&$ Mathematics with Applications}

\begin{document}

\begin{frontmatter}

\title{An improved immersed finite element particle-in-cell method for plasma simulation\tnoteref{mytitlenote}}
\tnotetext[mytitlenote]{This work was supported by National Natural Science Foundation of China [Grant numbers 10875034, 11175052];
Shenzhen Technology Project [Grant numbers JCYJ20150403161923511, JCYJ20150529115038093, JCYJ20160226201347750].}

\author[HITaddress]{Jinwei Bai}

\author[HITaddress]{Yong Cao}

\author[HITaddress,MSTaddress]{Yuchuan Chu\corref{mycorrespondingauthor}}
\cortext[mycorrespondingauthor]{Corresponding author}

\author[MSUaddress]{Xu Zhang}

\address[HITaddress]{Department of Mechanical Engineering $\&$ Automation, Harbin Institute of Technology, Shenzhen Graduate School, Shenzhen, Guangdong 518055, P. R. China}
\address[MSTaddress]{Department of Mathematics $\&$ Statistics, Missouri University of Science $\&$ Technology, Rolla, MO 65401, USA}
\address[MSUaddress]{Department of Mathematics $\&$ Statistics, Mississippi State University, Mississippi State, MS 39762, USA}

\begin{abstract}
The particle-in-cell (PIC) method has been widely used for plasma simulation, because
of its noise-reduction capability and moderate computational cost. The immersed finite element
(IFE) method is efficient for solving interface problems on Cartesian meshes, which is
desirable for PIC method. The combination of these two methods provides an effective tool
for plasma simulation with complex interface/boundary. This paper introduces an improved
IFE-PIC method that enhances the performance in both IFE and PIC aspects. For the electric
field solver, we adopt the newly developed partially penalized IFE method with enhanced
accuracy. For PIC implementation, we introduce a new interpolation technique to ensure the
conservation of the charge. Numerical examples are provided to demonstrate the features of
the improved IFE-PIC method.
\end{abstract}


\end{frontmatter}

\linenumbers





\section{Introduction}
There are two classes of methods for plasma simulation. The first one is the traditional dynamic simulation \cite{VlasovAA}, which
is mainly used to obtain the distribution function of particles change with time by solving the Vlasov equation. The second one is the particle simulation method \cite{JDawson,OCEldridge,1981HockneyEastwood}, which is used to track a large number of individual particles and to obtain the trajectory parameters and characteristics of plasma by statistical methods. Due to the enormous number of particles to be tracked and the limited computational resources, the development of the particle simulation method was quite slow. The particle simulation method entered into a rapid developing period \cite{PCBirch,SKimura,LASchwager} since Birdsall and Langdon \cite{LangdonAB_Birdsall} introduced the particle-in-cell (PIC) method which utilizes the finite-sized particle (or cloud) instead of a huge number of real particles. \vspace{2mm}

Immersed finite element (IFE) method is a finite element method for solving interface problems on uniform Cartesian meshes, which was first proposed by Li et al \cite{1998Li}. Different from classical finite element methods using body-fitted meshes, the mesh of IFE method is independent of the interface. However, the IFE basis functions around interfaces are modified to accommodate the interface jump conditions. The advantage of IFE method is that structured Cartesian meshes can be used to solve interface problems with arbitrary interface geometry. For problems with a moving interface, IFE methods are especially advantageous since there is no need to regenerated the solution meshes repeatedly
\cite{2013HeLinLinZhang,2013LinLinZhang2,2013LinLinZhang1}.
The IFE methods have been developed for solving the second-order elliptic equations \cite{2008HeLinLin, 2011HeLinLin, 2004LiLinLinRogers, 2003LiLinWu, 2015LinYangZhang1}, elasticity equations \cite{2010GongLi,2013LinSheenZhang,2012LinZhang}, Stokes equations \cite{2015AdjeridChaabaneLin}, to name only a few.\vspace{2mm}

In the past decade, IFE method has been used together with PIC method for plasma simulations \cite{YCao_YChu_XMHe_TLin,2017ChuHanCaoHeWang,2005KafafyLinLinWang,2005KafafyWangLin}. IFE method used as an electric field solver is performed on well-structured Cartesian meshes. This is particularly desirable for PIC method because tracking a large number of plasma particles can be efficiently done in a uniform structured mesh.  We refer to a few recent applications of IFE-PIC method for different types of particle simulations, such as ion thruster \cite{2015CaoChuWangCaoXiaZhang,2015JianChuCaoCaoHeXia}, hall thruster \cite{2015CaoLiShanCaoZheng}, and lunar surface environment \cite{2016HanWangHe,DHan_PWang}. Also, IFE-PIC method has been extended to handle unbounded interface problems with asymptotic boundary condition \cite{2011ChuCaoHeLuo} and periodic boundary condition \cite{2016CaoChuZhangZhang}. \vspace{2mm}


For the current IFE-PIC method, we noticed that there are two issues. First, the classical Galerkin IFE method is used as the field solver. As shown in \cite{2015LinLinZhang}, the classical IFE method is not accurate around the interface, because the IFE basis functions are discontinuous across the element boundaries, and the classical Galerkin formulation cannot control such discontinuity. Second, the particle interpolation method of PIC algorithm is imperfect. The conventional interpolation approach applied on the interface element often leads to the non conservation of charge, because it neglects the fact that some nodes of the interface cells are inside the conductors. Similar problems occur in applying the electric field force to the particles on the interface elements.\vspace{2mm}

In this paper, we introduce an improved IFE-PIC method that focuses on overcoming the problems mentioned above. As a remedy of discontinuity of IFE field solver, we adopt the newly developed partially penalized immersed finite element (PPIFE) method \cite{2015LinLinZhang} to improve the accuracy of IFE methods near interfaces. For PIC interpolation, we introduce a two-step approach for particle interpolation that preserves the charge conservation. Comparing with conventional charge distribution in PIC, we add a correction step that redistributes the quantity distributed to the nodes inside the conductor to the remaining nodes in order to maintain charge conservation. In addition, we use IFE basis functions to calculate the electric field and force on the interface elements. The new approach can calculate the motion of particles more accurately.\vspace{2mm}

The rest of the article is organized as follows. In Section 2, we recall the classical IFE method and PIC interpolation. In Section 3, we present our improved IFE-PIC method. The improvement in the IFE solver part is the PPIFE method with additional penalty terms. The improvement in the PIC part includes the new particle interpolation scheme and the new method for the force deposit. In Section 4, we present some numerical experiments to compare the performance of traditional IFE-PIC method and improved IFE-PIC method. Brief conclusions will be given in Section 5.

\section{Review of IFE-PIC Method}
In this section, we first recall the main steps in a typical IFE-PIC computational cycle. Then we will recall the classical IFE method and PIC interpolation that are widely used in the literature.

\subsection{Main Steps of IFE-PIC Method}
Real plasma particles are modeled as many macro-particles in the PIC method, and they follow the evolution of the orbits of individual particles in the self-consistent electromagnetic field. The field is then updated by solving the governing elliptic equation with discontinuous dielectric coefficients.
The IFE-PIC method is an iteration of solving for the electromagnetic field and particle motion until the steady state is achieved.\vspace{2mm}

In general, an IFE-PIC computational cycle consists of the following five steps:
\begin{description}
\item[Step 1. Initialization]
A series of initial settings of the simulation including domain, mesh, boundary condition, and initial position and velocity of particles must be set up.

\item[Step 2. Particle Push]
The motion of the particle is induced by the particles themselves and the applied external fields $\bm{E}$. The trajectory of an individual charged particle is obtained by integrating the Newton-Lorentz equation
\begin{equation}
  m\frac{d\textbf{v}}{dt}=q(\bm{E}+\textbf{v}\times\bm{B}),
\end{equation}
where $m$, $q$ and $\textbf{v}$ are the mass, the charge, and the velocity of the particle, respectively.  $\bm{B}$ denotes the static magnetic field.

\item[Step 3. Charge Deposit]
The change of the positions of particles leads to the change of the charge density  $\rho$ on each node. Thus, we need to calculate the charge density at each node according to the new positions of particles. The process of interpolating the particle charges on the discrete mesh points is called weighting.
In traditional PIC method, we use the ratio of the area of the rectangle formed by the opposite cell vertex and the particle to the area of element as the weighting.

\item[Step 4. Solving for Potential]
After obtaining the charging density $\rho$ on each node, the electric field should also be updated. The potential function $\Phi(\mathbf{x})$ can be described by the second-order Poisson's equation with discontinuous dielectric coefficient $\beta(\mathbf{x})$, which represents different types of material:
\begin{equation}
-\nabla\cdot(\beta\nabla\Phi)=\rho(\mathbf{x}).
\end{equation}
To solve this equation, we use IFE method as a field solver, because of its applicability of Cartesian mesh, which is desirable in the PIC simulation for fast tracking of particles' locations.

\item[Step 5. Force Deposit]
After solving for the potential $\Phi(\mathbf{x})$, we calculate the electric field $\bm{E}$ at each node by
\begin{equation}
  \bm{E}(\mathbf{x})=-\nabla\Phi(\mathbf{x}).
\end{equation}
Next, we need to deposit the electric field at nodes to the particles with arbitrary positions. The electric field of each node can be obtained using two-point difference method from the potential $\Phi(\mathbf{x})$ on each element node. The electric field of each particle can be obtained by interpolating.
\end{description}

We note that, in step 3, when the particle is located in an interface element, the current interpolation using the simple area-weighting will result in the non-conservation of charge. Also, in step 4, the classical IFE method may not be as accurate around interface as the rest of the domain due to the discontinuity of IFE basis functions. Finally, in step 5, the force deposit is not accurate in the interface element due to part of the element is the conductor. In this paper, our improvement of the IFE-PIC method will mainly focus on these three steps.

\subsection{Classical IFE Method for Interface Problems}
\begin{figure}[htb]
    \centering
    \includegraphics[angle=0, width=2.0in]{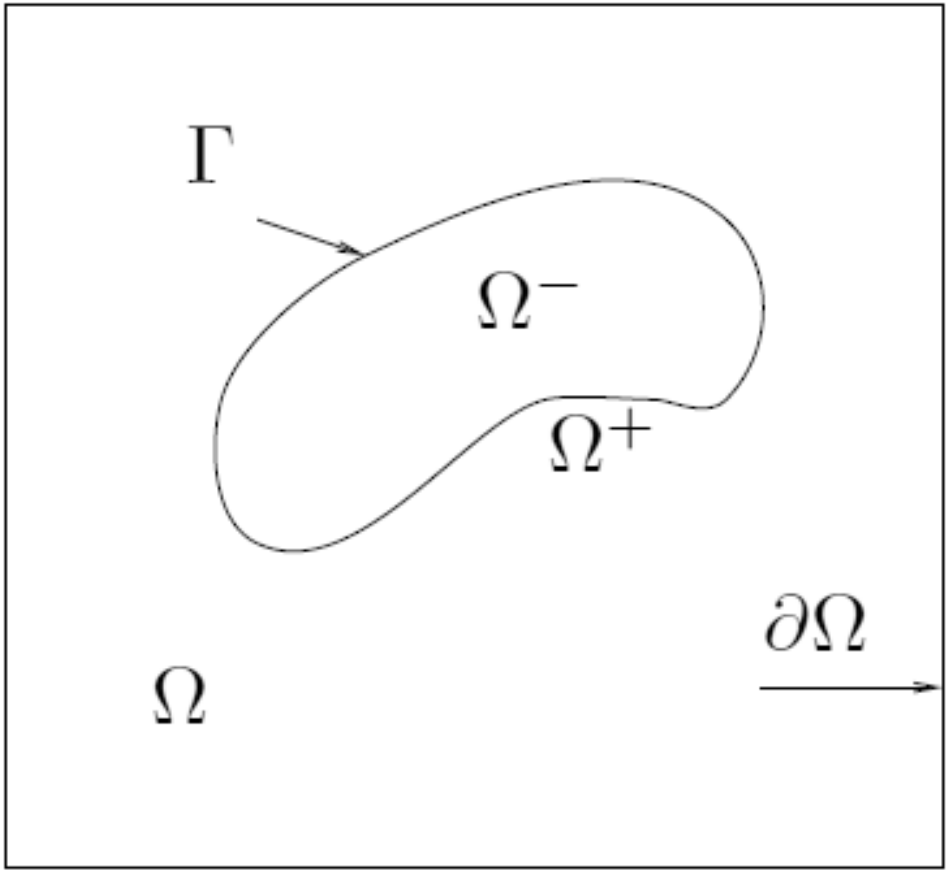}
    \caption{A sketch of the domain for the interface problem.}
    \label{Fig: domain}
\end{figure}

The electric field $\Phi(\mathbf{x})$ is assumed to be governed by the following second-order elliptic equation
\begin{eqnarray}
-\nabla\cdot(\beta\nabla\Phi)&=&f, ~~\text{in}~\Omega,\label{eq: pde}\\
\Phi&=&g,~~\text{on}~\partial\Omega.\label{eq: bc}
\end{eqnarray}
Here, we assume that $\Omega\subset \mathbb{R}^2$ is a rectangular domain separated by an interface curve $\Gamma$ into two sub-domains $\Omega^-$, $\Omega^+$ such that $\overline{\Omega}=\overline{\Omega^+\cup\Omega^-\cup\Gamma}$. See Figure \ref{Fig: domain} as an illustration.

The coefficient $\beta(\mathbf{x})$  is discontinuous across the interface. Without loss of generality, we assume $\beta$ is a piecewise constant function as follows
\begin{equation}\label{eq: coef}
  \beta(\mathbf{x})=\left\{
               \begin{array}{c}
                 \beta^-,~~\mathbf{x}\in\Omega^-, \\
                 \beta^+,~~\mathbf{x}\in\Omega^+,
               \end{array}
             \right.
\end{equation}
where $\mathbf{x} = [x,y]^t$.
Across the interface $\Gamma$, the following interface jump conditions are satisfied:
\begin{equation} \label{eq: jump}
  \jump{\Phi}_\Gamma=0,
\end{equation}
\begin{equation} \label{eq: flux jump}
  \bigjump{\beta\frac{\partial\Phi}{\partial n}}_\Gamma=0.
\end{equation}

\begin{figure}[!htp]
    \centering
    \includegraphics[angle=0, width=2.6in]{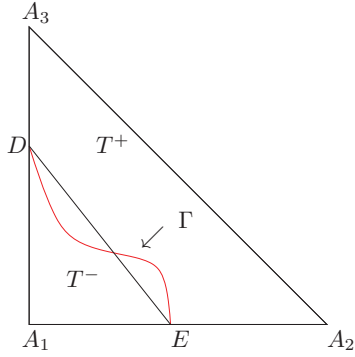}
    \caption{An interface triangle}
    \label{Fig: interface tri}
\end{figure}

Let $\mathcal{T}_h = \{T\}$ be a uniform triangular or rectangular mesh of the domain $\Omega$ with size $h$. If an element $T$ is cut through by the interface, it is called an interface element; otherwise, it is said to be a noninterface element. The sets of interface elements and noninterface elements are denoted by $\mathcal{T}_h^i$ and $\mathcal{T}_h^n$, respectively. Standard linear or bilinear finite element functions are used on all noninterface elements. Special piecewise-polynomial basis functions are constructed on interface element to accommodate the interface conditions. To be more precise, we use the linear IFE method as an example. Assume that $T=\triangle A_1A_2A_3$ is an interface triangle, and the interface curve $\Gamma$ intersects $T$ at two points, denoted by $D$ and $E$. The element $T$ is divided into two sub-elements $T^+$ and $T^-$ by the line segment $\overline{DE}$. See Figure \ref{Fig: interface tri} as an illustration. The local IFE basis functions $\phi_{i,T}$, $i=1,2,3$ on $T$ are defined as follows:
\begin{equation}\label{eq: ife basis}
  \phi_{i,T}(\mathbf{x})=\left\{
                 \begin{array}{lr}
                 \phi_{i,T}^+(\mathbf{x})=a_i^+x+b_i^+y+c_i^+,~~\mathbf{x}\in T^+,& \\
                 \phi_{i,T}^-(\mathbf{x})=a_i^-x+b_i^-y+c_i^-,~~\mathbf{x}\in T^-,&
                 \end{array}
              \right.
              \forall~i=1,2,3.
\end{equation}
They satisfy the following condition:
\begin{description}
\item[1. nodal-value conditions]
\begin{equation}\label{eq: nodal}
  \phi_{i,T}(A_j) = \delta_{ij} =
  \left\{
  \begin{array}{c}
  1,~~\text{if}~ i=j,\\
  0,~~\text{if}~ i\neq j,\\
  \end{array}
  \right.~~~\forall ~i,j = 1,2,3.
\end{equation}
\item[2. function-value continuity]
\begin{equation}\label{eq: cont}
    \phi_{i,T}^+(D)=\phi_{i,T}^-(D), ~~~
    \phi_{i,T}^+(E)=\phi_{i,T}^-(E), ~~~
  \forall ~i=1,2,3.
\end{equation}
\item[3. flux continuity]
\begin{equation}\label{eq: flux}
  \beta^+\frac{\partial\phi_{i,T}^+}{\partial \textbf{n}}=\beta^-\frac{\partial\phi_{i,T}^-}{\partial \textbf{n}},  ~~\forall ~i=1,2,3,
\end{equation}
where $\textbf{n}$ is the normal vector of $\overline{DE}$.
\end{description}

The local linear IFE space is defined as
\begin{equation}
\tilde{\mathbb{P}}_1(T) = span\{\phi_{1,T},\phi_{2,T},\phi_{3,T}\}.
\end{equation}
Denote the set of interior nodes on $\mathcal{T}_h$ by $\mathcal{N}_h$. On each node $\mathbf{x}_i$, $i=1,2,\cdots, N$, we define the global linear IFE basis function $\Phi_i$ such that
\begin{equation}
\Phi_i(\mathbf{x}_j) = \delta_{ij}, ~~~\forall~ i,j = 1,2,\cdots, N,
\end{equation}
and
\begin{equation}
\Phi_i|_{T} \in \tilde{\mathbb{P}}_1(T),~~\forall~ T\in \mathcal{T}^i_h,~~~\text{and}~~~
\Phi_i|_{T} \in {\mathbb{P}}_1(T) = span\{1,x,y\},~~\forall~ T\in \mathcal{T}^n_h.
\end{equation}
The global linear IFE space is formed as $S_h = span\{\Phi_i: i=1,2,\cdots, N\}$. \\

The classical (Galerkin) IFE method is to find $u\in S_h$ such that
\begin{equation}
a(u_h,v_h) = (f,v_h),~~~\forall~v_h\in S_h,
\end{equation}
where
\begin{equation}
a(u,v) = \sum_{T\in\mathcal{T}_h} \int_T \beta\nabla u\cdot\nabla v d\mathbf{x},~~~~
(f,v) = \int_\Omega f v d\mathbf{x}.
\end{equation}
The construction of the bilinear IFE spaces on rectangular meshes are similar, and we refer to \cite{2008HeLinLin, 2011HeLinLin, 2012LinZhang} for more details.

\subsection{Particle Interpolation in Traditional PIC Method}
\begin{figure}[!htp]
    \centering
    \subfigure[Non interface element]{
    \includegraphics[angle=0, width=2.2in]{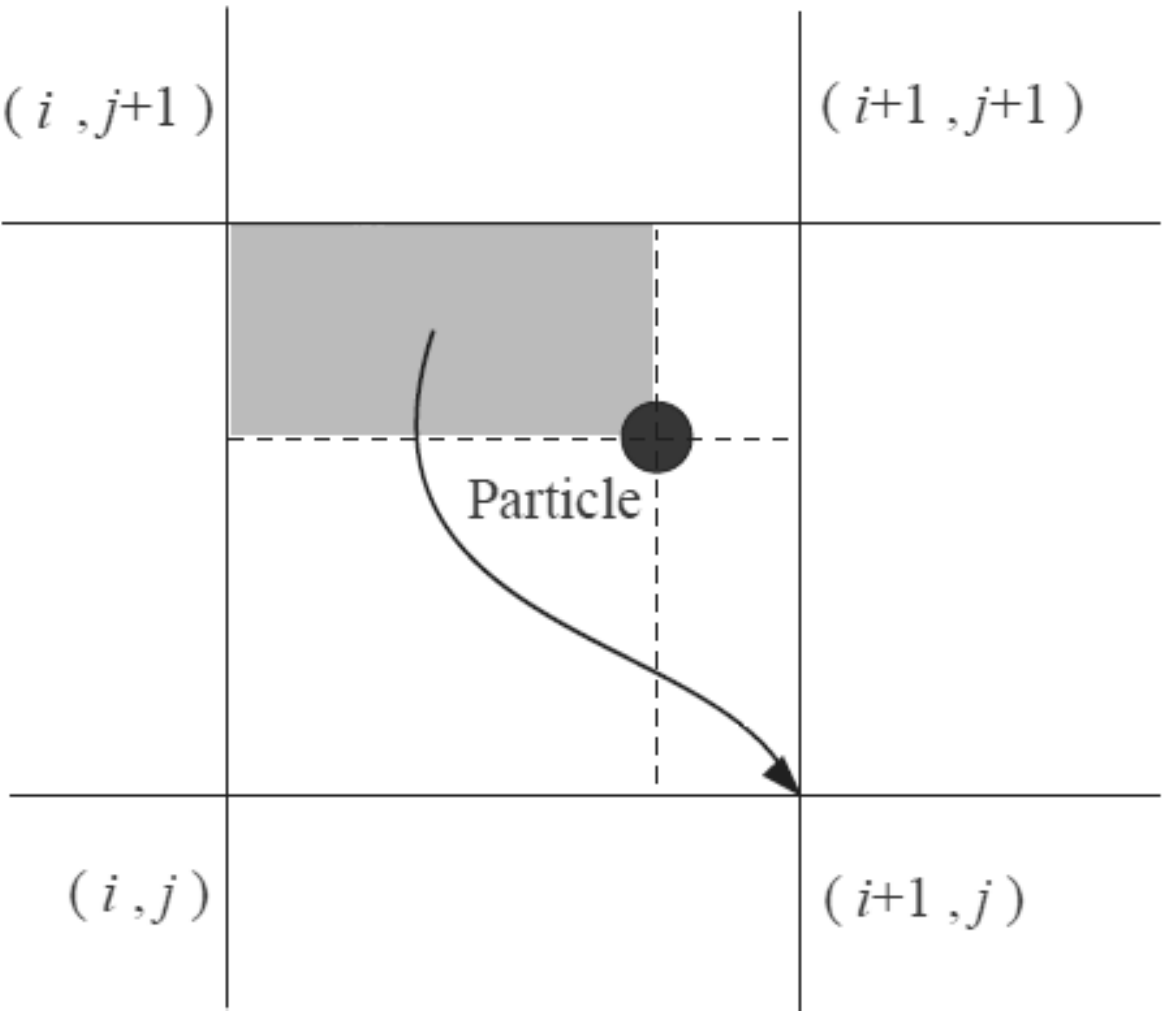}
    \label{Fig:Non-interface element}
    }~~~~~
    \subfigure[Interface element]{
    \includegraphics[angle=0, width=2.2in]{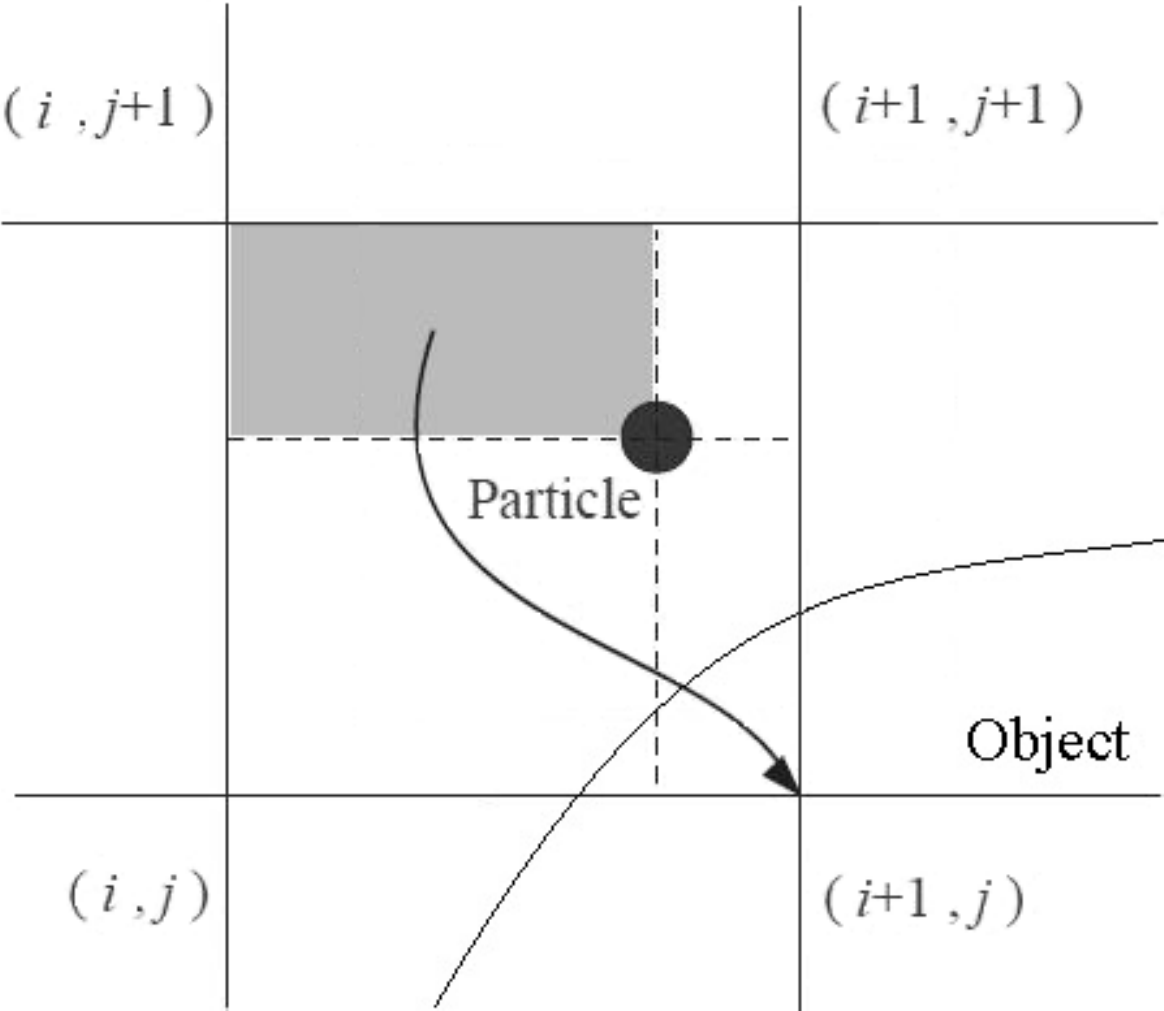}
    \label{Fig:Interface element}
    }
    \caption{Standard PIC deposition scheme of particle charge in a 2D simulation domain.}
    \label{Fig: interpolation.}
\end{figure}

In the PIC method, particle interpolation is required at Step 3 and Step 5 described in Section 2.1.
First, in charge deposit in Step 3,  as shown in Figure \ref{Fig:Non-interface element}, the portion of the total particle charge assigned to a certain cell vertex is proportional to the area of the rectangle formed by the opposite cell vertex and the particle.
Thus, the interpolation of the particle $P$, located at position $X_P$, to the node $X_{i,j}$ can be calculated by
\begin{equation} \label{eq: interpolation}
  q_{i,j}=q_P\frac{S{(X_{i+1,j+1},X_P)}}{S{(X_{i,j},X_{i+1,j+1})}},
\end{equation}
where $q_{i,j}$ is the amount of charge of node $(i,j)$, $q_P$ is the amount of charge of particle $P$, and $S(A,B)$ is the area of the rectangle whose diagonal is $\overline{AB}$.

On non-interface elements, the interpolation \eqref{eq: interpolation} works well. However, on interface elements, it will face some difficulty, as shown in Figure \ref{Fig:Interface element}, the charge should not be assigned to the node $X_{i+1,j}$, which is inside the conducting object. Directly applying \eqref{eq: interpolation} on interface elements will cause the non-conservation for the total charge and charge density.\vspace{2mm}

Moreover, in Step 5 the electric field is required to be interpolated at the particles' positions, so that the electric field force can be obtained. The potential $\Phi$ at grid points are solved by IFE method with the appropriate boundary conditions. Then, the electric field $\bm{E} = [E_x,E_y]^t$ can be obtained from the potential $\Phi$ using the following equation
\begin{equation}\label{eq: potential}
  \bm{E}(x,y)=-\nabla\Phi(x,y).
\end{equation}
The conventional approach for this interpolation uses the finite difference form for $E_x$ and $E_y$:

\begin{equation}\label{eq: linear interpolation}
(E_x)_{i,j} = \frac{\Phi_{i-1,j}-\Phi_{i+1,j}}{2\Delta x},~~
(E_y)_{i,j} = \frac{\Phi_{i,j-1}-\Phi_{i,j+1}}{2\Delta y},
\end{equation}

where $\Delta x$ and $\Delta y$ are mesh sizes in the $x$- and $y$- directions, respectively. After the electric field is obtained, the forces caused by the field at mesh nodes can be deposited to the arbitrary particles positions.

When a particle is located in the interface element, the calculation \eqref{eq: linear interpolation} is apparently inaccurate, because of the electric field of the four points that located at the different locations of the object are discontinuous.

\section{Improved IFE-PIC Method}
In this section, we present an improved IFE-PIC method. The improvement involves both the IFE solver and the PIC particle interpolation.

\subsection{Partially Penalized IFE Method}
Due to the discontinuity of IFE basis functions, classical IFE method using Galerkin formulation may generate large errors around the interface. Partial penalized immersed finite element (PPIFE) method is introduced in \cite{2015LinLinZhang} that is known to greatly improve the accuracy of IFE solution around interface. The main idea of this method is to add penalty terms on interface edges to reduce the negative impact of the discontinuity introduced by IFE functions. In our improved IFE-PIC method, we adopt this new PPIFE scheme as the new field solver.

To present the PPIFE method, we need to introduce a few notations. Let $\mathcal{E}_h$ be the set of all interior edges of the mesh $\mathcal{T}_h$. If an edge $e\in\mathcal{E}_h$ intersects with the interface curve $\Gamma$, we call it an interface edge; otherwise a noninterface edge. The sets of interface edges and noninterface edges are denoted by $\mathcal{E}_h^i$ and $\mathcal{E}_h^n$, respectively. For each interior edge $e\in \mathcal{E}_h$, it must be shared by two adjacent elements, denoted by $T_{e,1}$ and $T_{e,2}$. For a function $u$ defined on $T_{e,1}\cup T_{e,2}$, the average and jump of $u$ on $e$ are defined as follows
\begin{equation}\label{eq: averjump}
  \aver{u}_e=\frac{1}{2}\Big((u|_{T_{e,1}})_e+(u|_{T_{e,2}})_e\Big),~~~~
  \jump{u}_e=(u|_{T_{e,1}})_e-(u|_{T_{e,2}})_e.
\end{equation}

The PPIFE method for solving \eqref{eq: pde} - \eqref{eq: flux jump} is to find $u_h\in S_h$ such that
\begin{equation}\label{eq: ppife}
   a_\epsilon (u_h,v_h) = (f, v_h),~~~\forall v_h\in S_h,
\end{equation}
where
\begin{eqnarray}
a_\epsilon(u,v) &=&  \sum_{T\in\mathcal{T}_h} \int_T \beta\nabla u\cdot\nabla v d\mathbf{x}
 - \sum_{e\in\mathcal{E}_h} \int_e \aver{\beta \nabla u\cdot \mathbf{n}}\jump{v}ds \nonumber\\
&& + \epsilon\sum_{e\in\mathcal{E}_h} \int_e \aver{\beta \nabla v\cdot \mathbf{n}}\jump{u}ds +
 \sum_{e\in\mathcal{E}_h} \int_e \frac{\sigma_e}{|e|}\jump{u}\jump{v}ds. \label{eq: ppife bilinear form}
\end{eqnarray}
Here, $\epsilon$ can be chosen as $-1$, $0$, or $1$, which corresponds to symmetric, incomplete, or non-symmetric PPIFE methods, respectively.  $\sigma_e$ is a positive penalty parameter. As shown in \cite{2015LinLinZhang}, for non-symmetric PPIFE method, the scheme \eqref{eq: ppife} is stable as long as $\sigma_e >0$ on every edge $e$. For symmetric and incomplete PPIFE methods, \eqref{eq: ppife} is stable when $\sigma_e $ is large enough.

\subsection{Charge-Conservative Particle Interpolation}
As shown in Section 2, the particle interpolation technique \eqref{eq: interpolation} is not accurate for particles in the interface elements. Also, the total charge
is not conservative due to part of the nodes are inside the objects such as metal or ceramic materials. In this subsection, we introduce a two-step interpolation that can maintain the charge conservation.
\subsubsection{Improved Algorithm for Charge Deposit}

\begin{figure}[!htp]
    \centering
    \subfigure[One node inside the object]{
    \includegraphics[angle=0, width=1.9in]{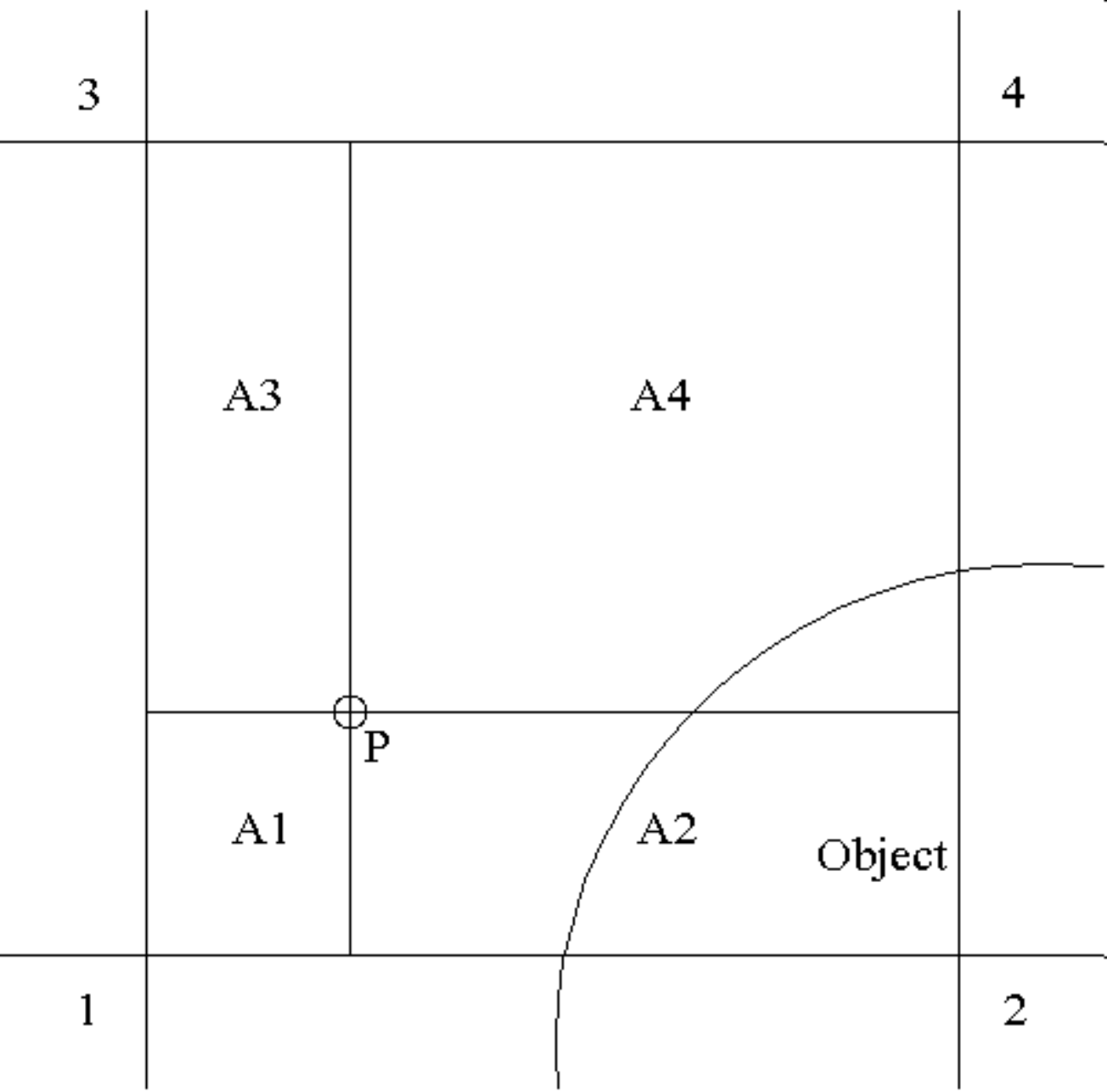}~
    \label{Fig:One node inside the object}
    }
    \subfigure[Two nodes inside the object]{
    \includegraphics[angle=0, width=1.9in]{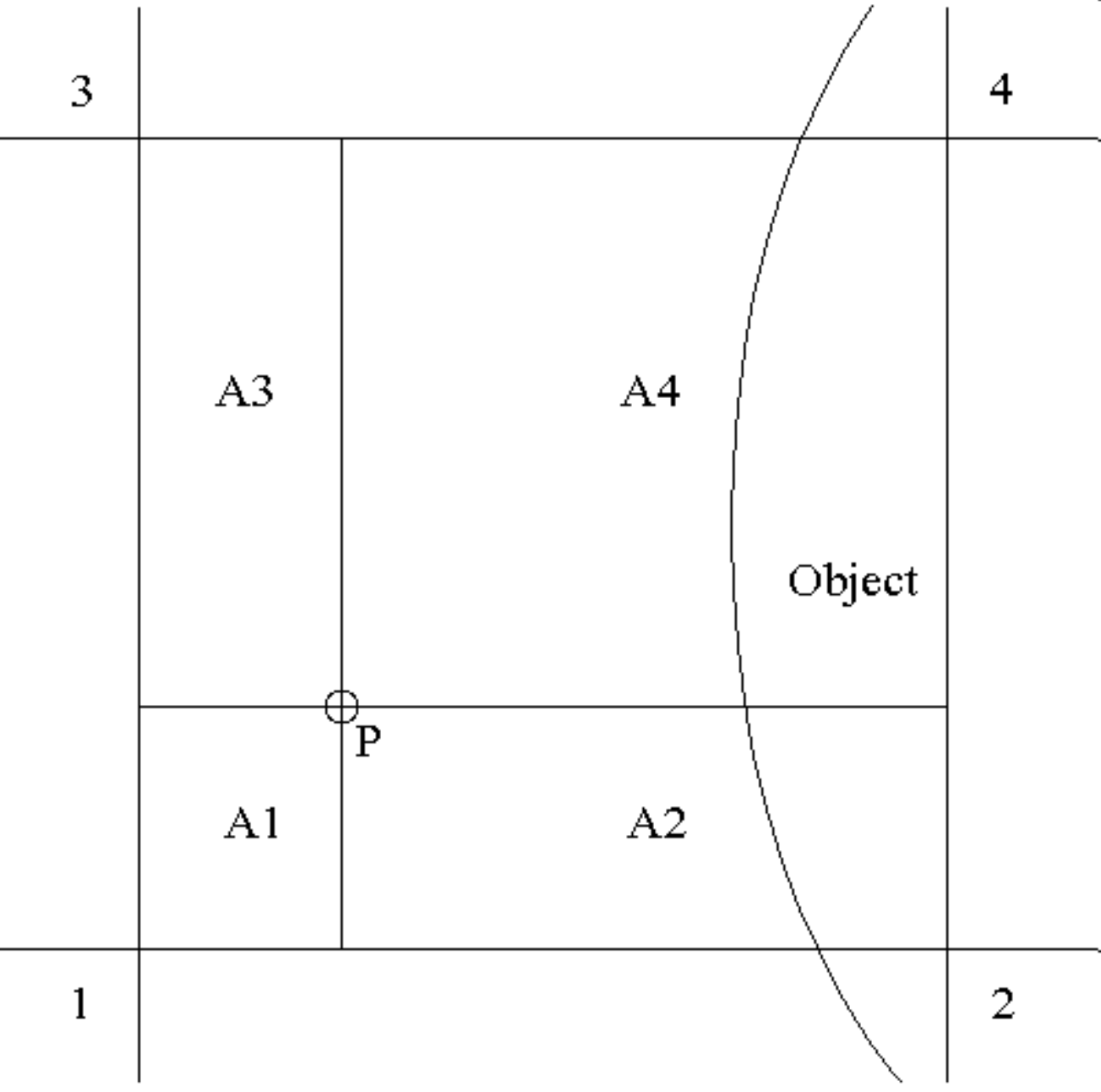}~
    \label{Fig:Two nodes inside the object}
    }
    \subfigure[Three nodes inside the object]{
    \includegraphics[angle=0, width=1.9in]{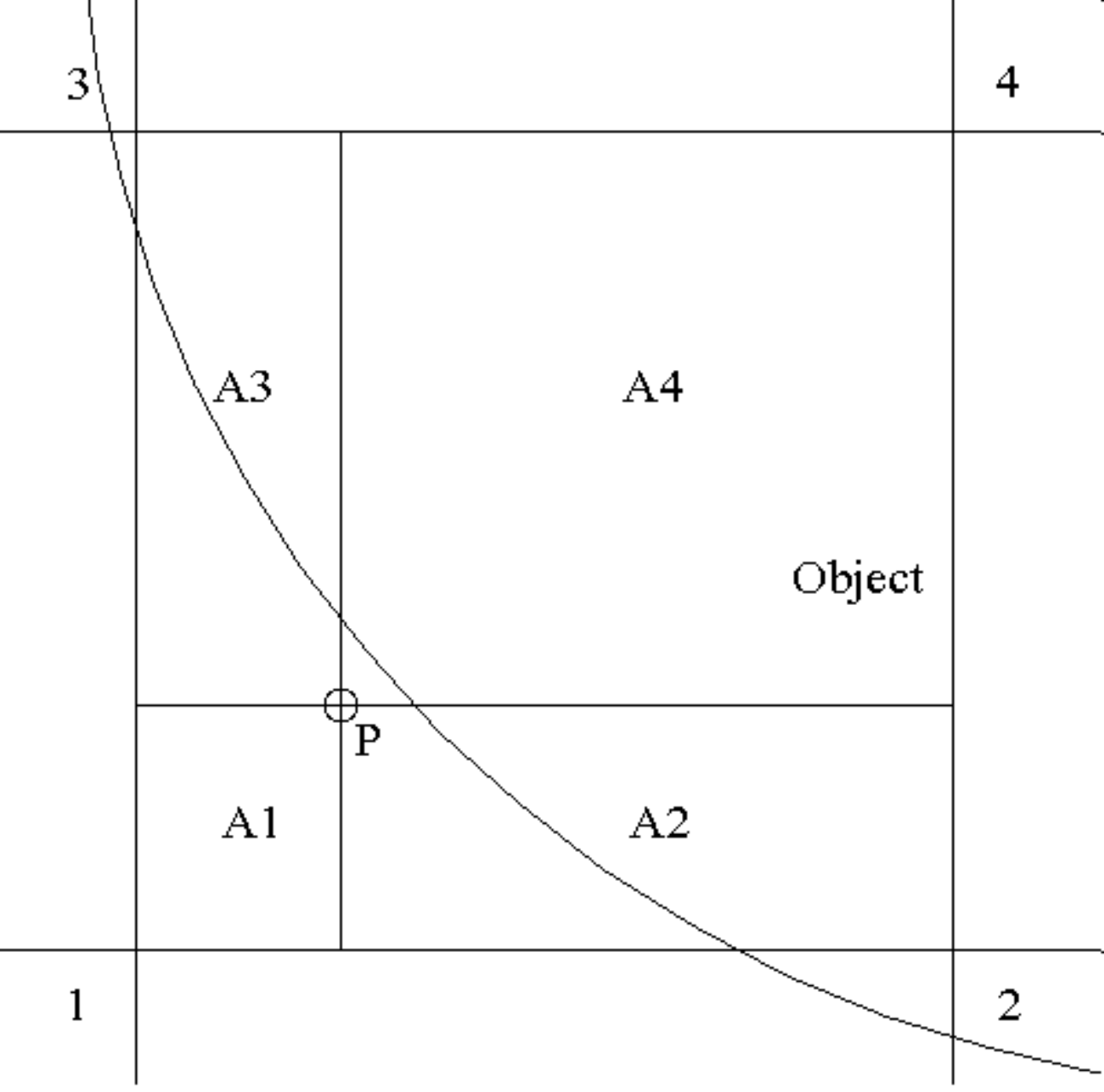}
    \label{Fig:Three nodes inside the object}
    }
    \caption{New deposit algorithm of particle charge.}
    \label{fig: interpolation}
\end{figure}

We note that if a particle is in a noninterface element, it is suitable to use the standard charge distribution \eqref{eq: interpolation} with area weight method. If a particle $P$ is in an interface element, we introduce a correction step to redistribute the charge assigned to nodes inside the conductor to other nodes. This step ensures the conservation of the charge in the simulation domain. The complete two-step interpolation is as follows.
\begin{description}
\item[Step 1: Initial distribution] Interpolate the charge $q_{_P}$ of the particle to four nodes of the element using the standard area-weight approach:
\begin{equation}\label{eq: step1}
  q_i^{(1)}=w_i q_{_P},~~~i = 1,2,3,4,
\end{equation}
where
\begin{equation*}
  w_1 = \frac{A_4}{A},~~w_2 = \frac{A_3}{A},~~w_3 = \frac{A_2}{A},~~w_4 = \frac{A_1}{A}.
\end{equation*}
Here $A_i$, $i=1,2,3,4$ are the areas of small rectangles, illustrated in Figure \ref{fig: interpolation}.
\item[Step 2: Re-distribution] Correct the charge distribution to nodes inside the object.
\begin{itemize}
\item \textbf{Case 1}:  There is only one node (e.g. node $\#$2) in the object, as shown in Figure \ref{fig: interpolation}(a). The charge initially distributed to node $\#$2 should be redistributed to the other three nodes. The correction procedure is the following:
\begin{eqnarray*}
q_1^{(2)} &=&q_1^{(1)}+\frac{w_1}{w_1+w_3+w_4} q_2^{(1)}, \\
q_2^{(2)} &=&0, \\
q_3^{(2)} &=&q_3^{(1)}+\frac{w_3}{w_1+w_3+w_4} q_2^{(1)}, \\
q_4^{(2)} &=&q_4^{(1)}+\frac{w_4}{w_1+w_3+w_4} q_2^{(1)}.
\end{eqnarray*}

\item \textbf{Case 2}: If two nodes are inside the object, as shown in Figure \ref{fig: interpolation}(b), then the correction procedure becomes
\begin{eqnarray*}
q_1^{(2)} &=&q_1^{(1)}+\frac{w_1}{w_1+w_3} (q_2^{(1)}+q_4^{(1)}), \\
q_2^{(2)} &=&0, \\
q_3^{(2)} &=&q_3^{(1)}+\frac{w_3}{w_1+w_3} (q_2^{(1)}+q_4^{(1)}), \\
q_4^{(2)} &=&0.
\end{eqnarray*}

\item \textbf{Case 3}: If three nodes are inside the object, as shown in Figure \ref{fig: interpolation}(c), then
\begin{eqnarray*}
q_1^{(2)} &=&q_1^{(1)}+ q_2^{(1)}+q_3^{(1)}+q_4^{(1)}, \\
q_2^{(2)} &=&0, \\
q_3^{(2)} &=&0, \\
q_4^{(2)} &=&0.
\end{eqnarray*}

\end{itemize}
\end{description}

We note that if we use the traditional PIC charge interpolation (without correction) on interface elements, the total charge of nodes is less than the charge of particle, that is
\begin{equation}
\sum_{i\in \Lambda}{q_i^{(1)}}\neq q_{_P},
\end{equation}
where $\Lambda$ denotes the indices of nodes outside the object.
However, adding the correction, the interpolation satisfies that
\begin{equation}
\sum_{i\in \Lambda}{q_i^{(2)}}= q_{_P}.
\end{equation}
This clearly shows that the new interpolation scheme preserves the charge-conservation. \vspace{2mm}

Moreover, this new interpolation scheme is robust regarding the location of the particle. That is, if the particle touches the boundary of  the interface element, the algorithm is still valid.

\subsubsection{Improved Algorithm for Force Deposit}

The potential distribution at grid points is obtained by solving linear algebraic equations contained in the immersed finite element method, and then for the electric potential at point $(x,y)$ in an element $T$ can be approximated using the following equation:

\begin{equation}
  \Phi_T(x,y)=\sum_{i=1}^{d}u_i\phi_{i,T}(x,y).
\end{equation}
Here, $u_i$ are the numerical solutions of the potential at the vertices of $T$, and $\phi_{i,T}$ are local IFE basis functions. The degree of freedom $d=3$ for linear IFE method, and $d=4$ for bilinear IFE method.

Suppose the electrical field is in $\Omega^+$, and the conductor is in $\Omega^-$. Then the electric field $\textbf{E}$ at a particle $P(x,y)$ can be obtained by
\begin{equation}
      \begin{array}{ccc}
      E_x(x,y)&=&-\left(u_1\dfrac{\partial \phi_1^+(x,y)}{\partial x}+u_2\dfrac{\partial \phi_2^+(x,y)}{\partial x}+u_d\dfrac{\partial \phi_d^+(x,y)}{\partial x}\right), \vspace{1mm}\\
      E_y(x,y)&=&-\left(u_1\dfrac{\partial \phi_1^+(x,y)}{\partial y}+u_2\dfrac{\partial \phi_2^+(x,y)}{\partial y}+u_d\dfrac{\partial \phi_d^+(x,y)}{\partial y}\right).
      \end{array}
\end{equation}

\begin{figure}[htb!]
    \centering
    \includegraphics[angle=0, width=.8\textwidth]{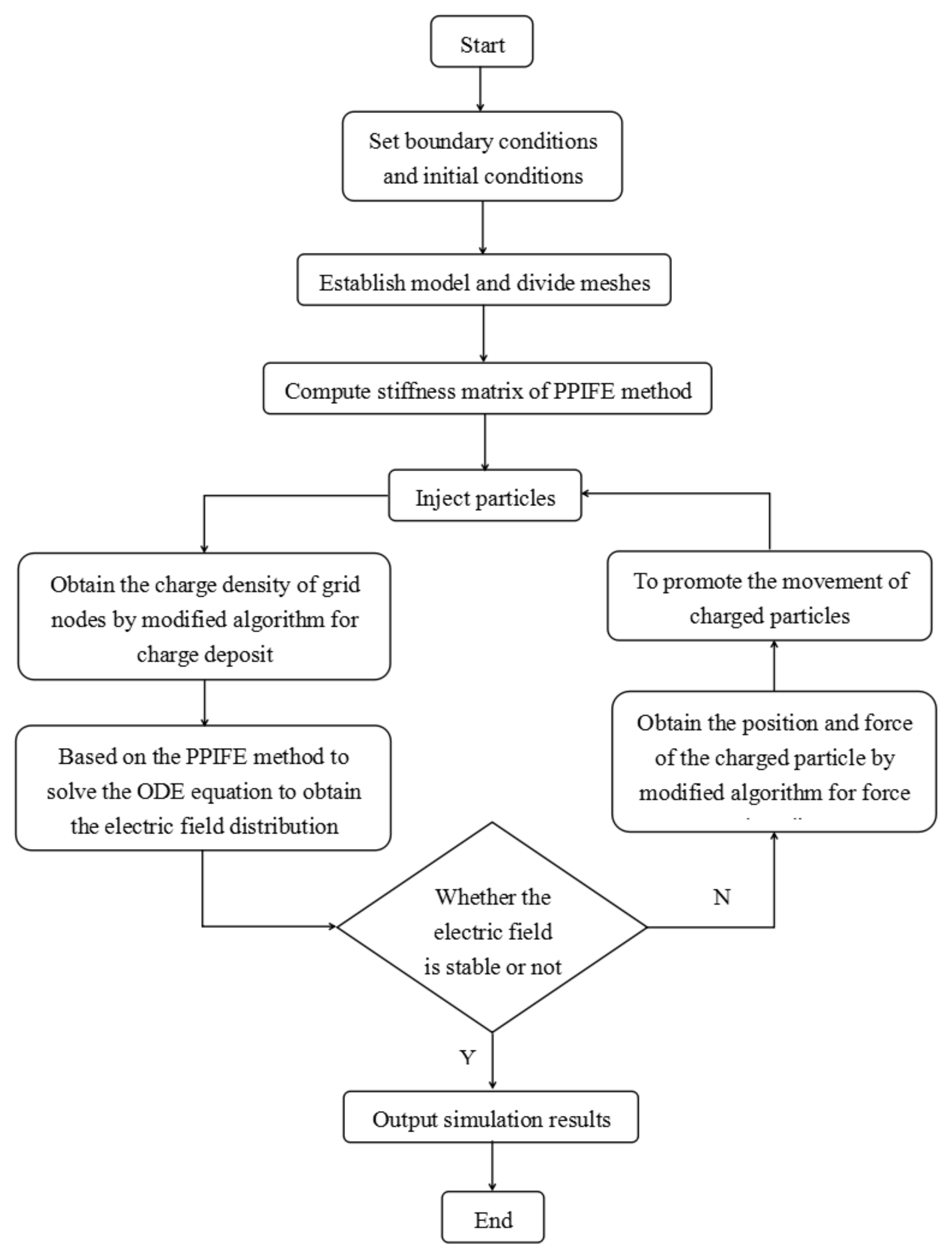}
    \caption{Flow chart of the improved IFE-PIC method.}
    \label{fig: workflow}
\end{figure}

For points located in $T^-$, the electric field can be calculated in the same way. This method is more accurate than linear interpolation when the internal and external potential of objects are discontinuous. With the new PPIFE field solver, and the improved PIC interpolation technique on particles and force, the workflow of our improved IFE-PIC algorithm can be summarized in Figure \ref{fig: workflow}.

\section{Numerical Examples}
In this section, we present some numerical examples to demonstrate the features of the improved IFE-PIC method.

\begin{figure}[!htp]
\centering
\includegraphics[width=3.0in]{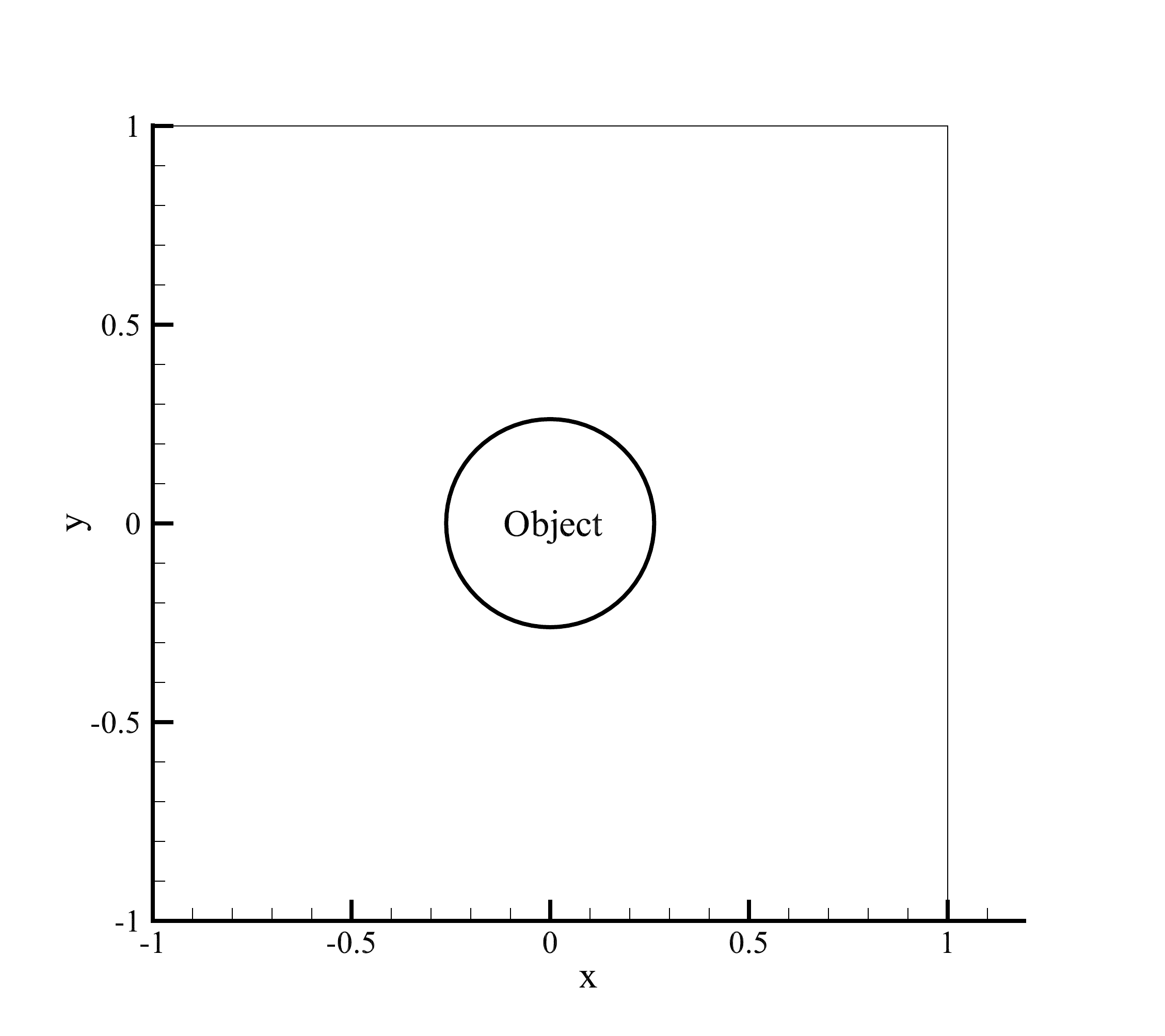}~
\includegraphics[width=3.0in]{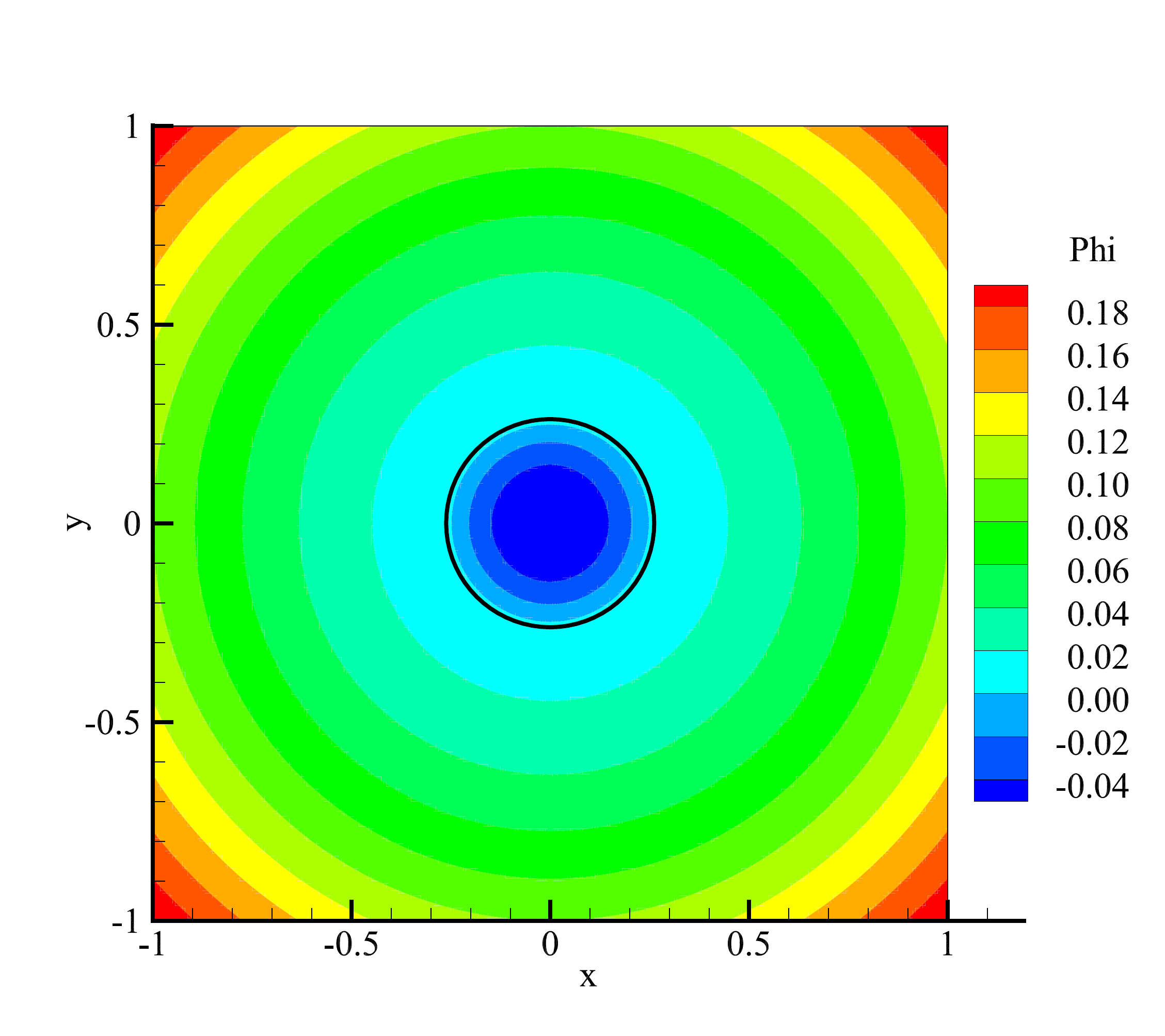}
\caption{Simulation domain $\Omega_s$ (left) and the analytical solution (right)}
\label{fig:domain}
\end{figure}

We set up a test problem of conducting cylinder with background plasma particles.
Let the two-dimensional simulation domain be $\Omega_s=[-1,1]\times[-1,1]$. The center and the radius of the cylinder are set to be $O = [0,0]$ and $r_0 =\pi/12$, respectively, as shown in Figure \ref{fig:domain}. For plasma particles, we load $1635841~ (=1279^2)$ uniformly distributed particles $P_{i,j}$ into the simulation domain and the positions of particles are given by
\begin{equation}
P_{i,j} =  [-1+\frac{i}{1280}, -1+\frac{j}{1280}],~~~i,j =1,2,...,1279.
\end{equation}

Assume that the loaded particles are electrons, and the charge density on the mesh point is $\rho=-4$, so we can get the charge of each particle. Then we remove all the particles inside the cylinder. We choose appropriate Dirichlet boundary conditions and source function so that the analytical solution of this simulation problem is
\begin{equation}\label{eq:exact}
  \Phi(x,y)= \left\{
  \begin{array}{ll}
  \Phi^-(x,y) =  \dfrac{1}{\beta^-} r^2+(\dfrac{1}{\beta^+}-\dfrac{1}{\beta^-})r_0^2,~~&\text{if}~~r\leq r_0,\vspace{1mm}\\
  \Phi^+(x,y)=\dfrac{1}{\beta^+}r^2,&\text{if}~~r>r_0.
  \end{array}\right.
\end{equation}
where $r = \sqrt{x^2+y^2}$, and $[\beta^-,\beta^+]=[1,10]$. The analytical solution \eqref{eq:exact} of the potential in the simulation domain is shown in Figure \ref{fig:domain}.

\subsection{Comparison of PIC Interpolations}
\begin{figure}[htp]
    \centering
    \includegraphics[angle=0, width=2.84in]{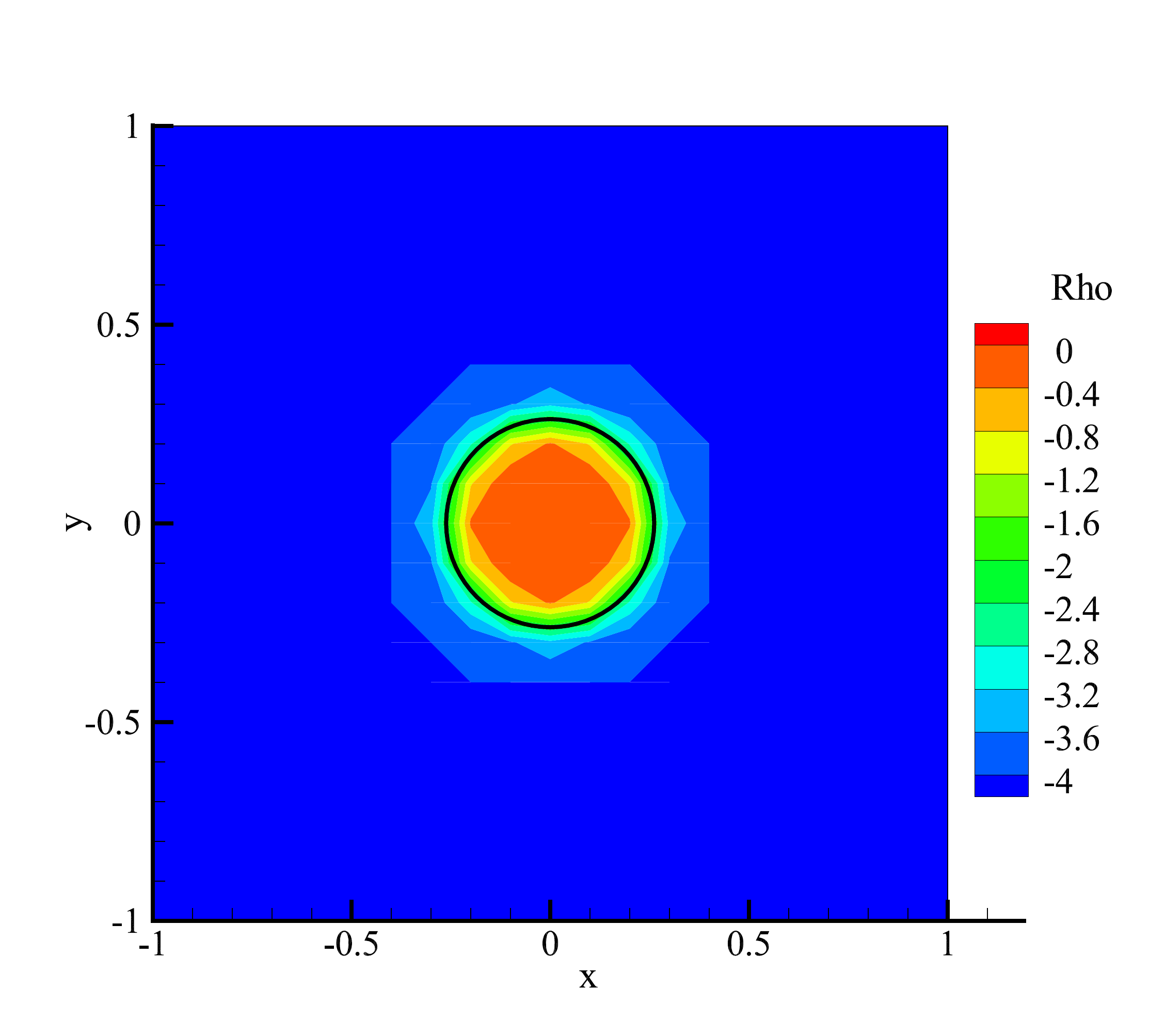}~~
    \includegraphics[angle=0, width=2.84in]{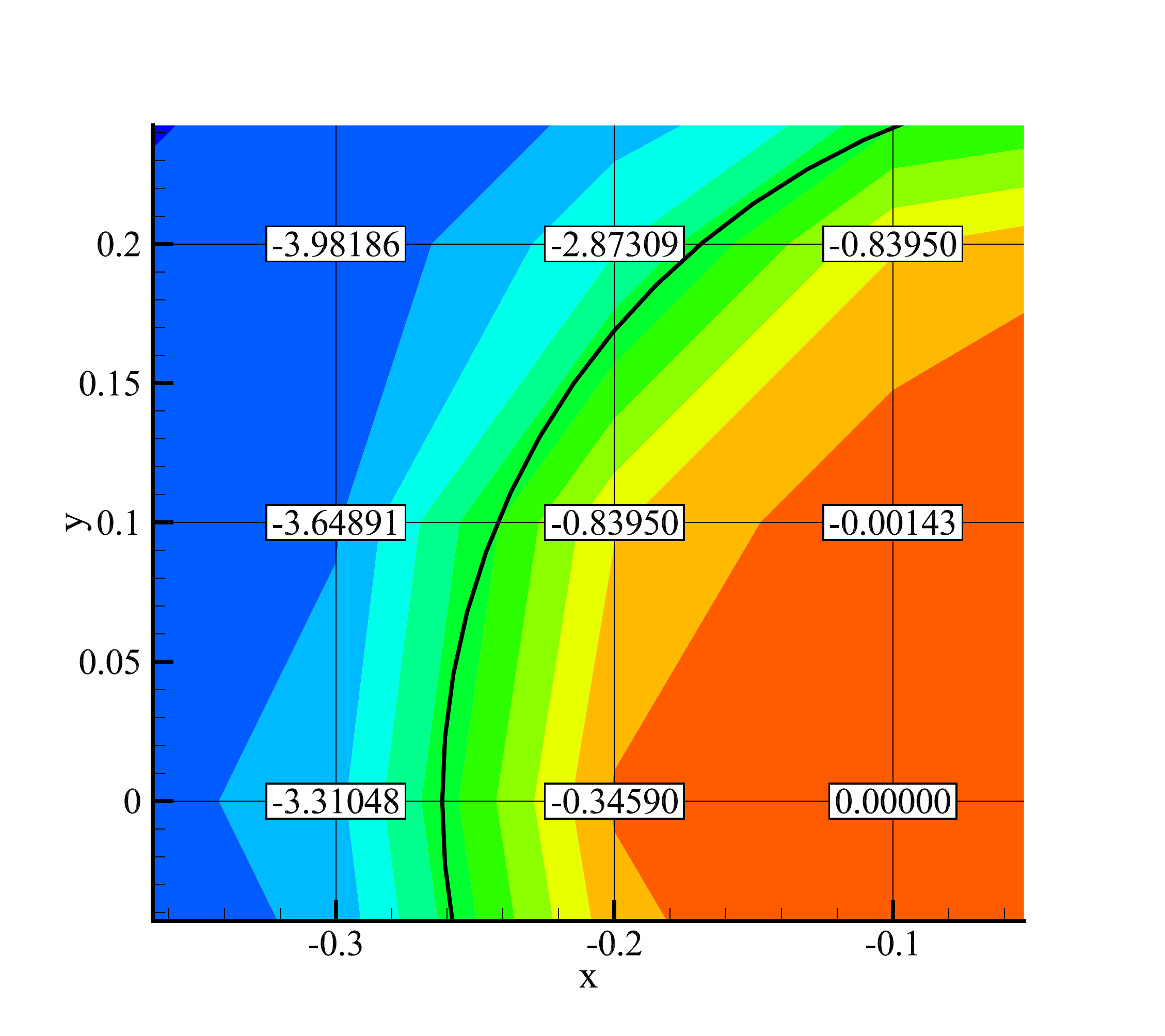}
        \caption{Density distribution using standard PIC interpolation}
    \label{fig: interpolation standard}
\end{figure}

We first compare the performance of traditional PIC charge interpolation with the new interpolation method. The charge densities on mesh points can be obtained by depositing the physical quantities (charge of electrons) from the particle locations onto mesh nodes. Figure \ref{fig: interpolation standard} shows the charge density distribution with standard charge deposit algorithm on a $20\times20$ Cartesian mesh. It can be shown that all the charge of interface element is less than the prescribed value $\rho =-4$, because the charge of particles are deposited to several mesh points inside the object.

\begin{figure}[htp]
    \centering
     \includegraphics[angle=0, width=2.84in]{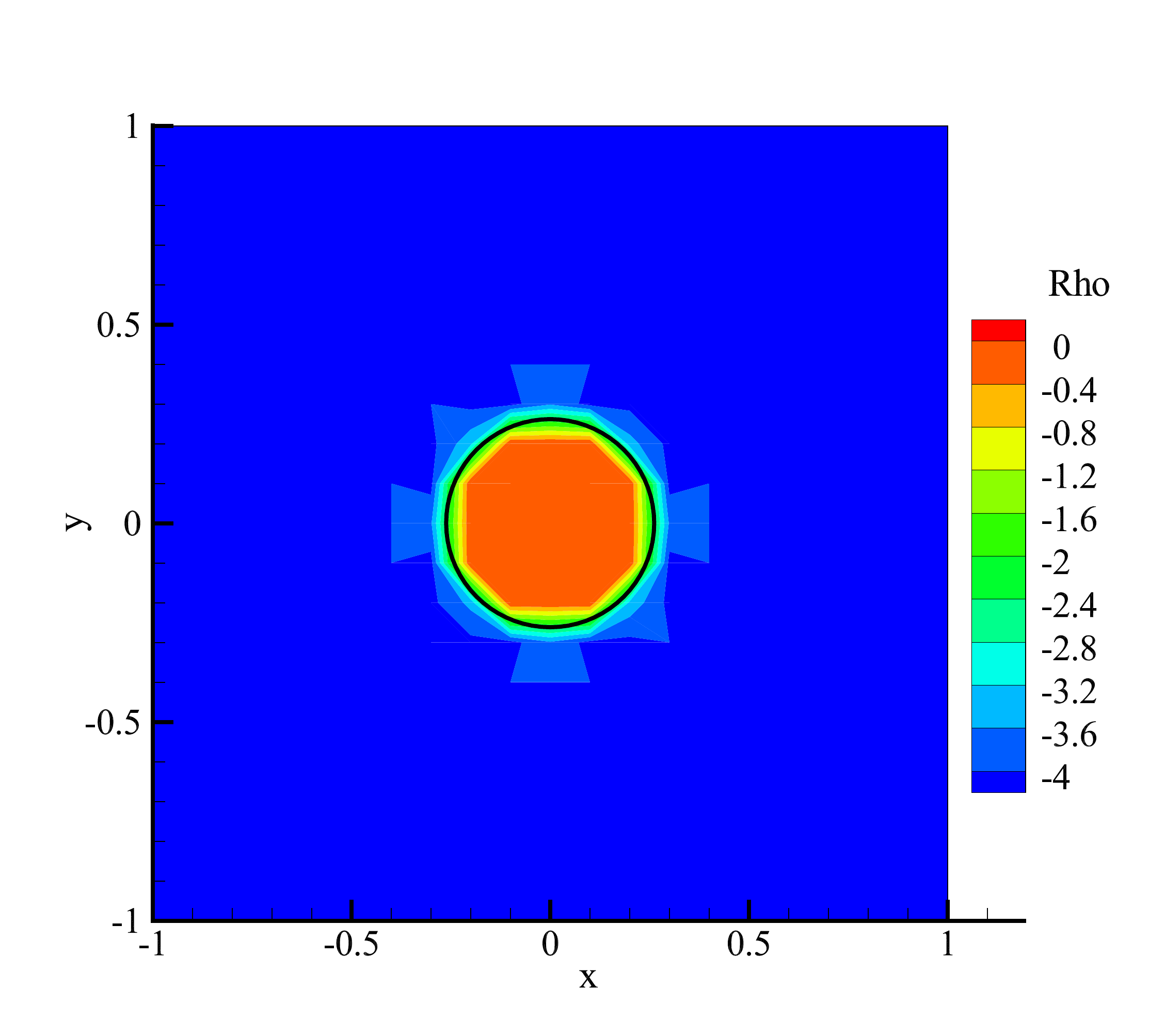}~~
     \includegraphics[angle=0, width=2.84in]{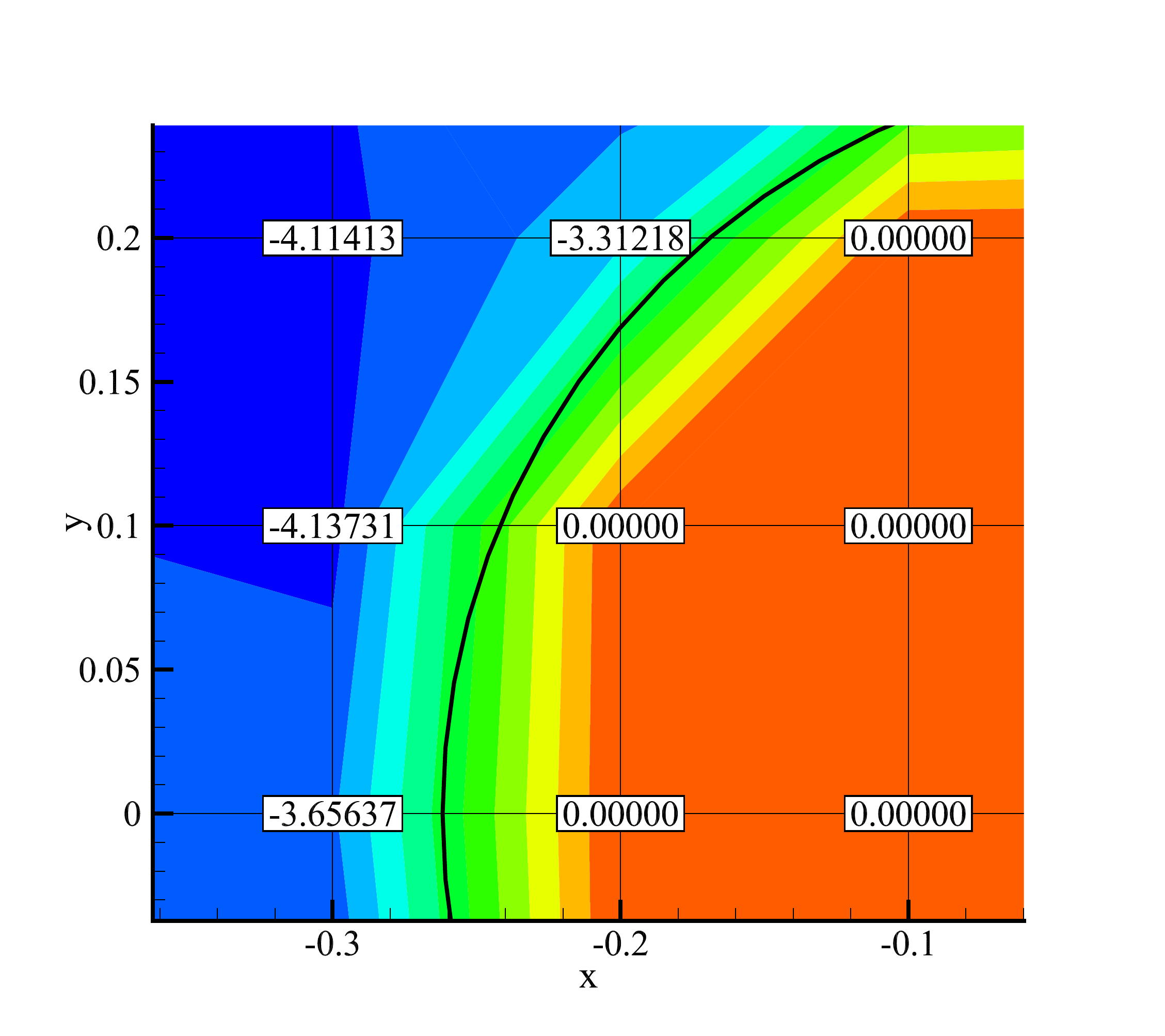}
     \caption{Density distribution using improved PIC interpolation}
    \label{fig: interpolation improved}
\end{figure}

Using the improved PIC algorithm introduced in Section 3.2, the charge density distribution in simulation domain $\Omega_s$ is shown in Figure \ref{fig: interpolation improved}. On the right side of Figure \ref{fig: interpolation improved}, we show a zoom-in plot around interface. It can be shown that the charge on the nodes of interface elements are not unanimously less than $-4$. On some nodes (e.g. [-0.3,0.1], and [-0.3,0.2]), the change quantities are actually greater than $-4$. This is because of the redistribution step we added, so that the charge of particles will not be deposited to any mesh points inside the conductor. Moreover, this new interpolation algorithm ensures the conservation of the charge in the calculation domain.

\subsection{Comparison of IFE-PIC Performance}
\begin{figure}[htb!]
    \centering
    \includegraphics[angle=0, width=2.84in]{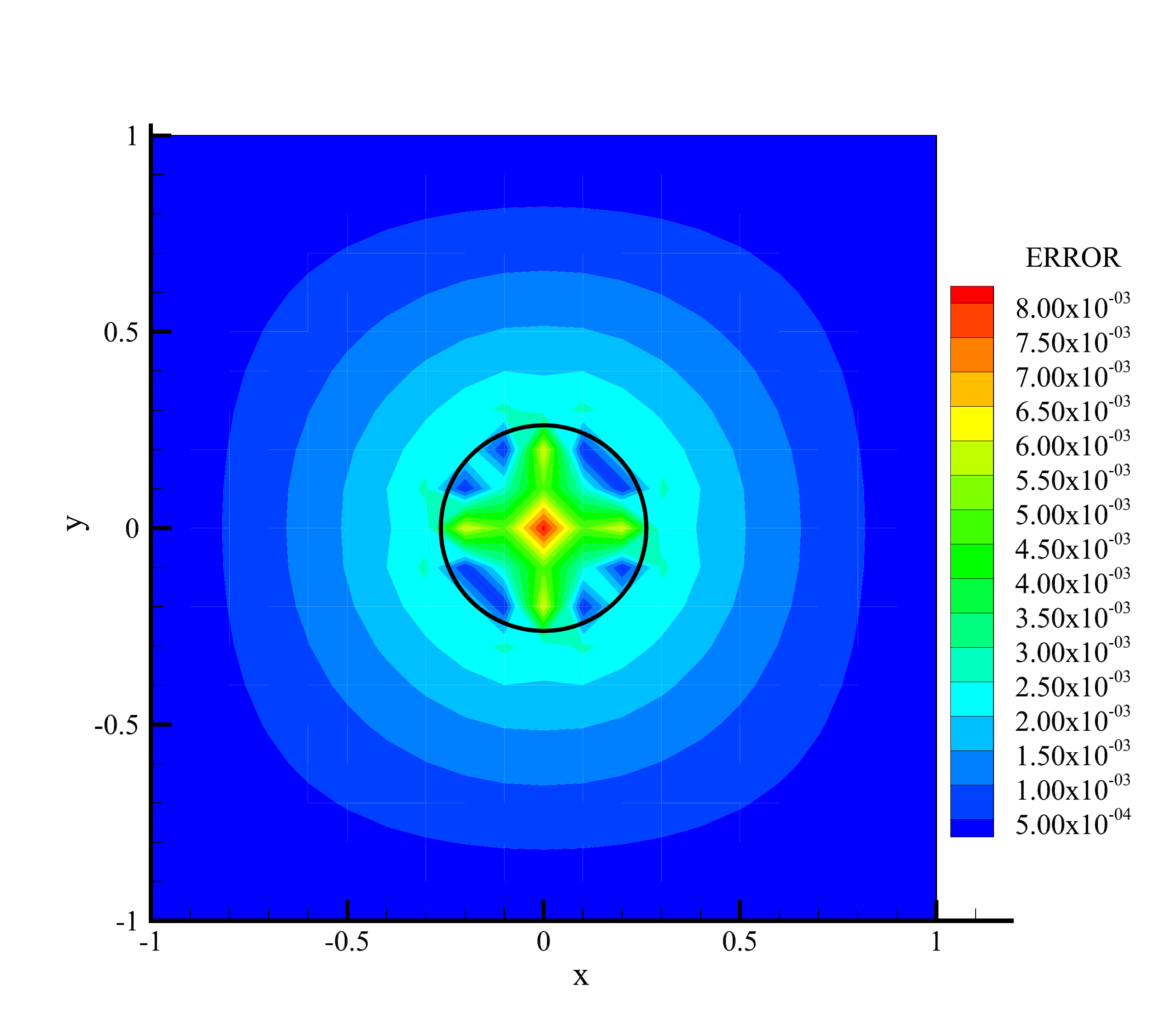}~~
        \includegraphics[angle=0, width=2.84in]{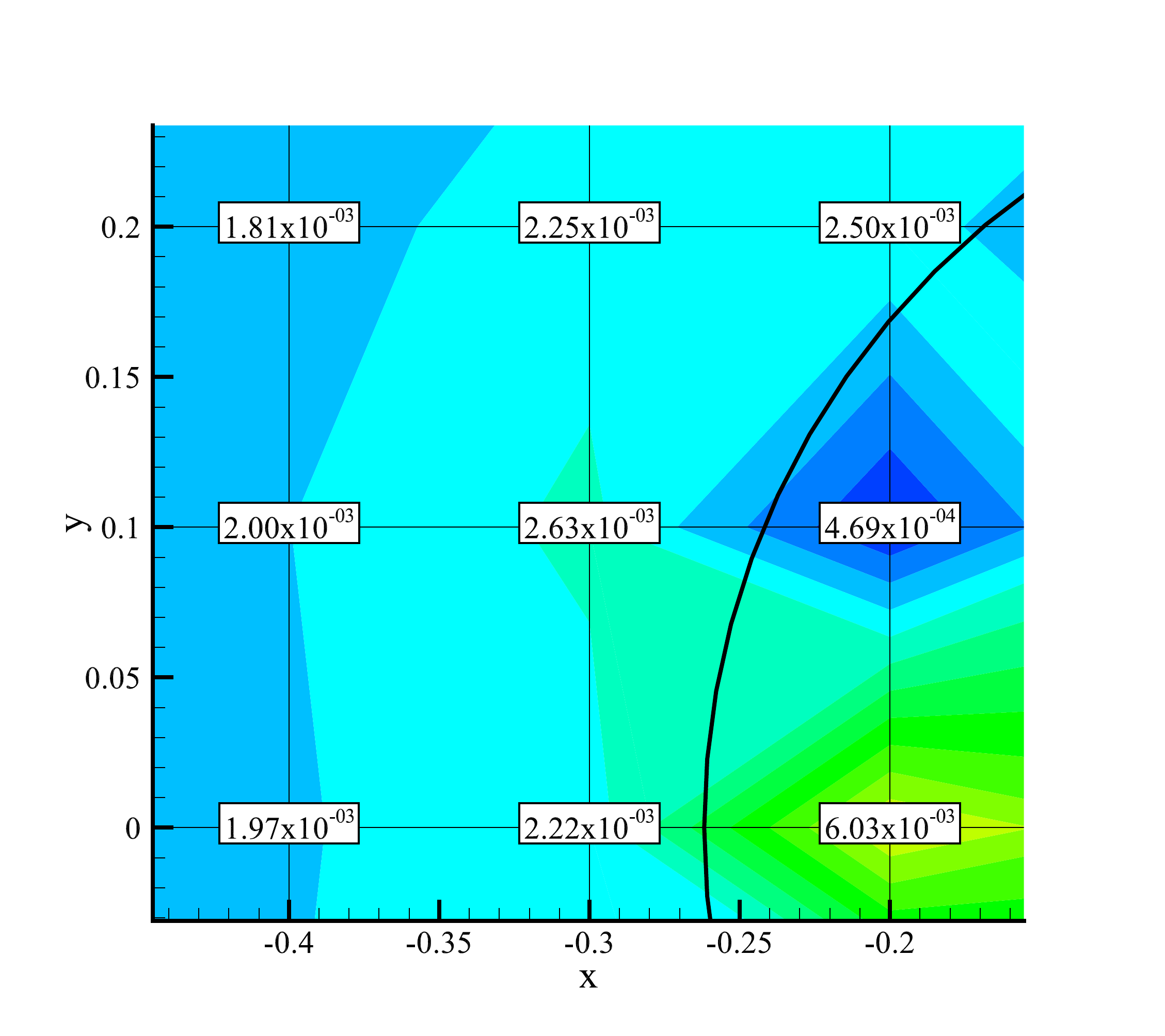}
    \caption{Error of traditional IFE-PIC solution on a Cartesian mesh of $\Omega_s$.}
    \label{fig: error standard}
\end{figure}

\begin{figure}[htb!]
    \centering
    \includegraphics[angle=0, width=2.84in]{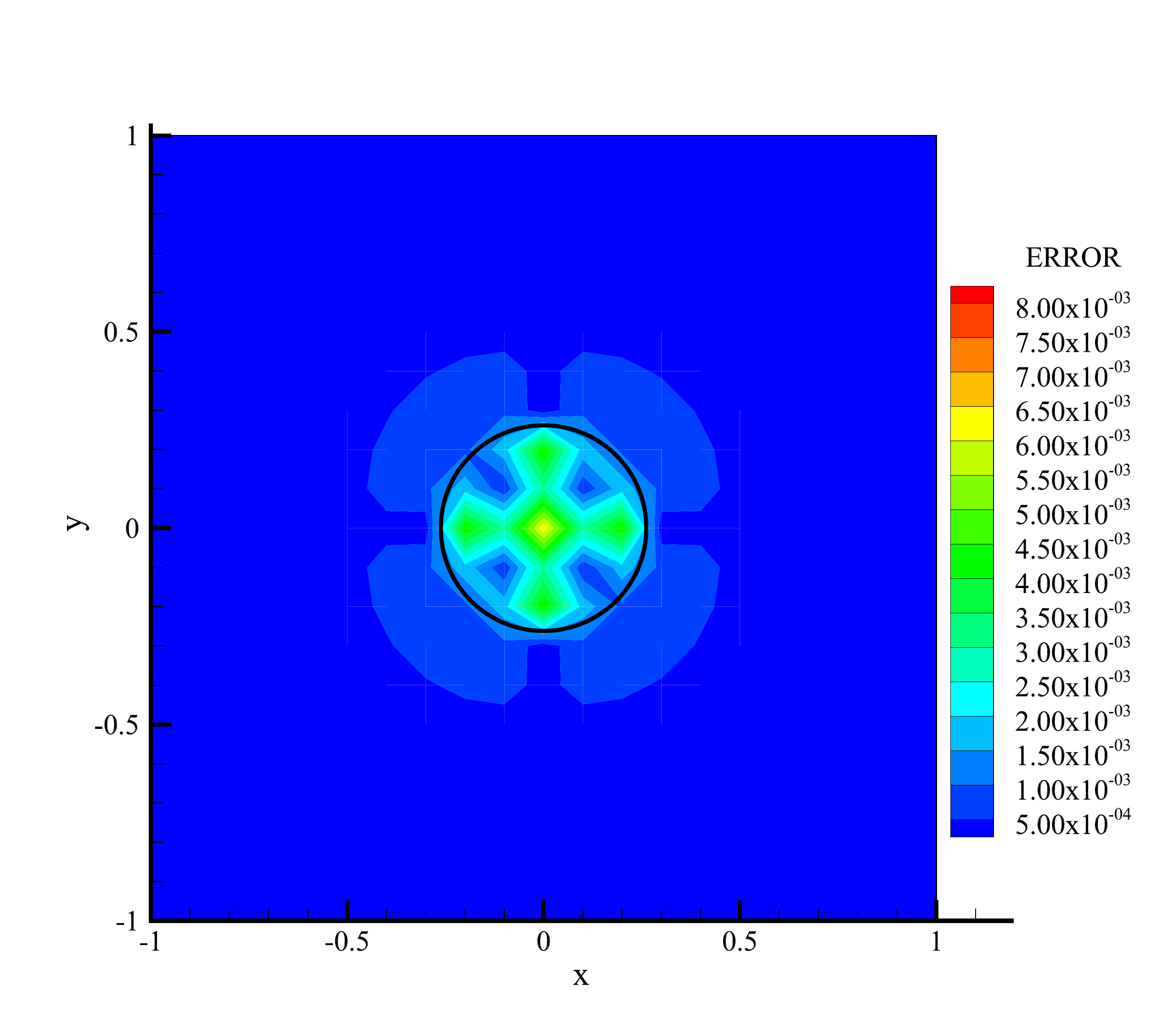}~~~
        \includegraphics[angle=0, width=2.84in]{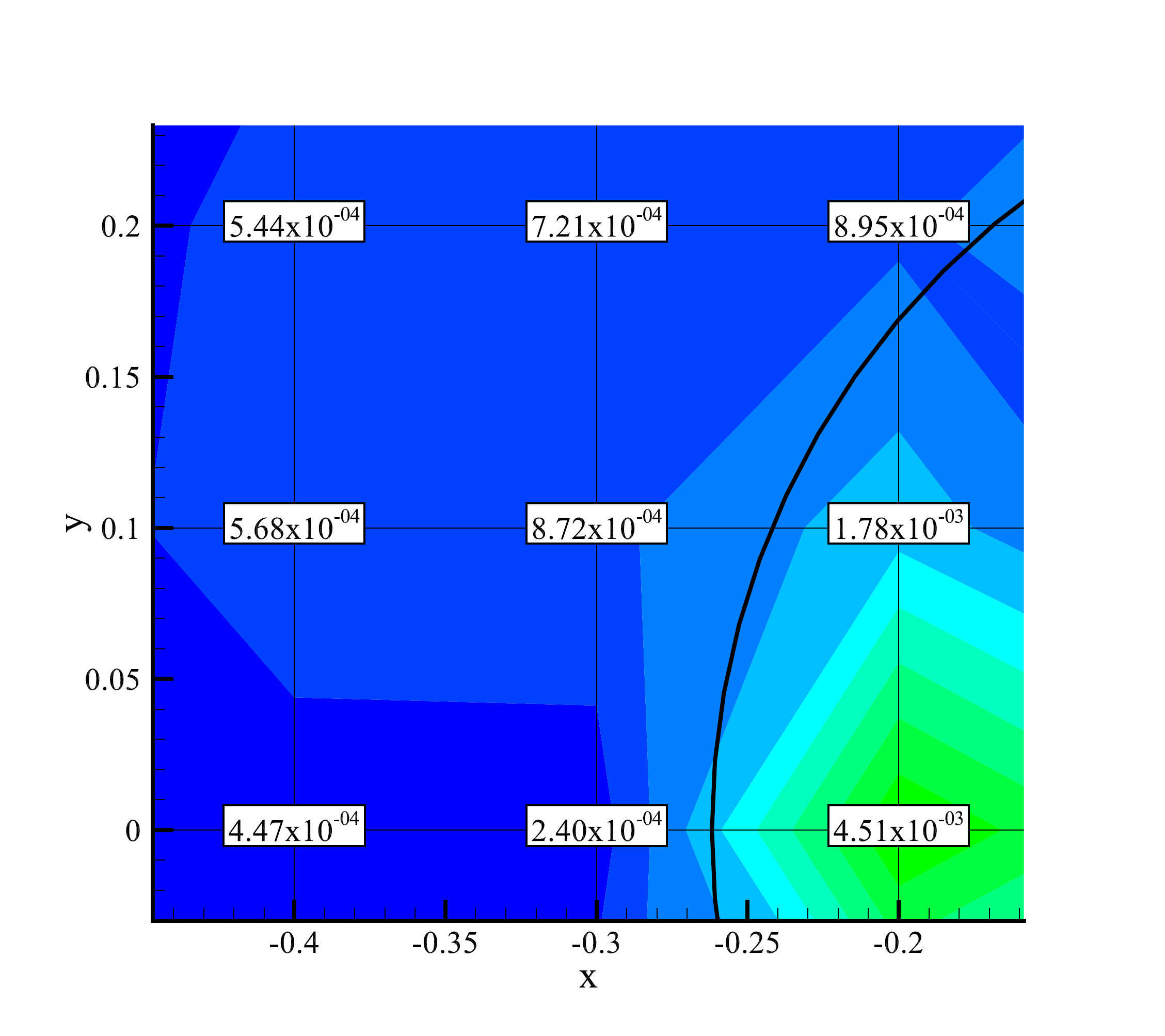}
            \caption{Error of improved IFE-PIC solution on a Cartesian mesh of $\Omega_s$.}
    \label{fig: error improved}
\end{figure}

Next, we combine the new PIC algorithm with the new electrical field solver, PPIFE scheme. We compare the performance of traditional and improved IFE-PIC
schemes. In Figure \ref{fig: error standard} and Figure \ref{fig: error improved}, we plot the numerical errors of electrical field by traditional and improved IFE-PIC methods, respectively. It can be easy observed that the new IFE-PIC method has much smaller error than traditional IFE-PIC method especially around the interface.

\subsection{Sensitivity of the Number of  Particles}
In this test, we focus on the sensitivity of our IFE-PIC method to the number of particles. Note that
on non-interface elements, regardless of the number of particles within each element, the final charge of interpolation to nodes are same, i.e., $\rho=-4$. However, on interface elements, the element is split into two parts. The particles located in the conductor will not contribute to the charge distribution, as a result, the number of particles in an interface element should have greater influence on the charge distribution results.

\begin{figure}[!htp]
    \centering
    \includegraphics[angle=0, width=1.6in]{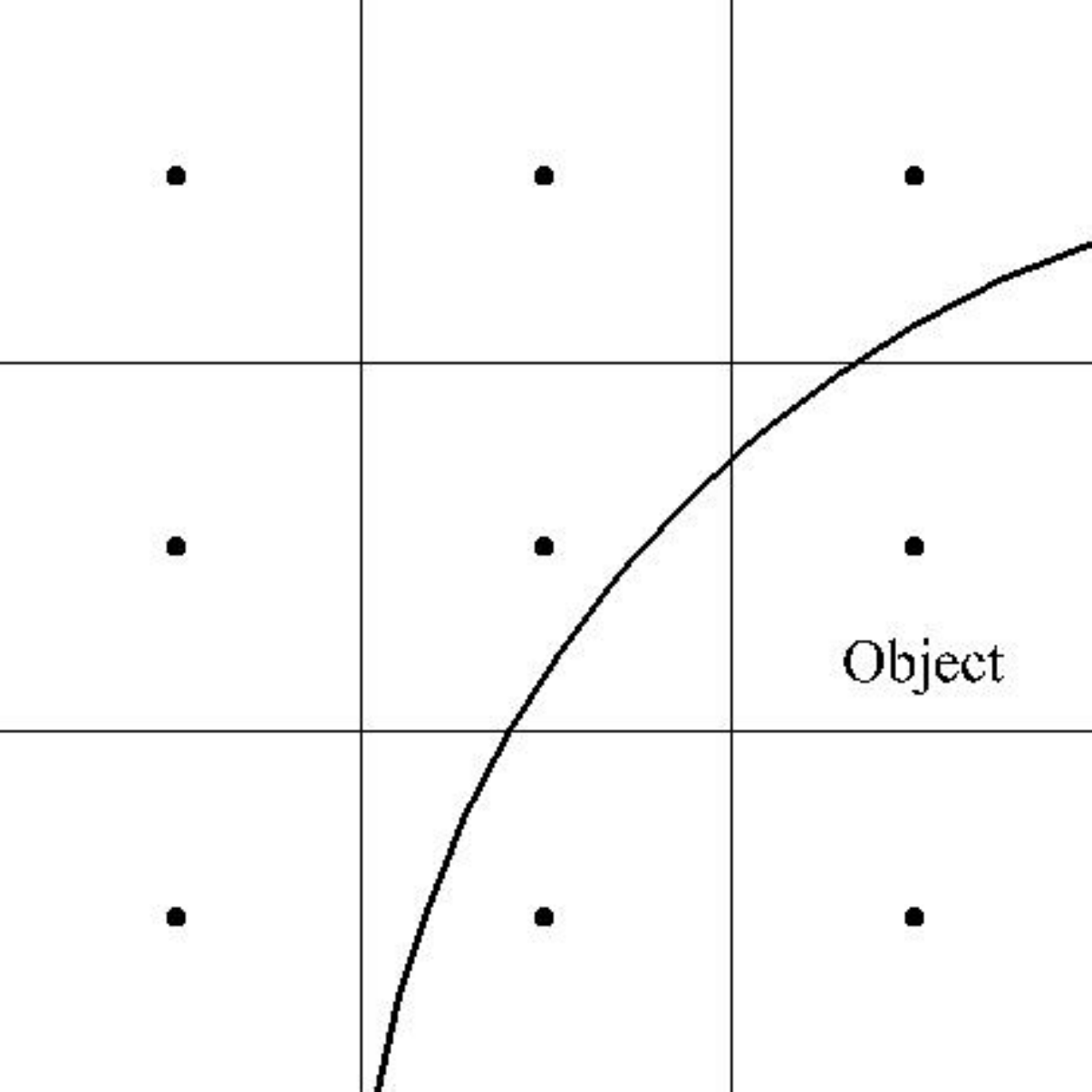} ~~~~~~
    \includegraphics[angle=0, width=1.6in]{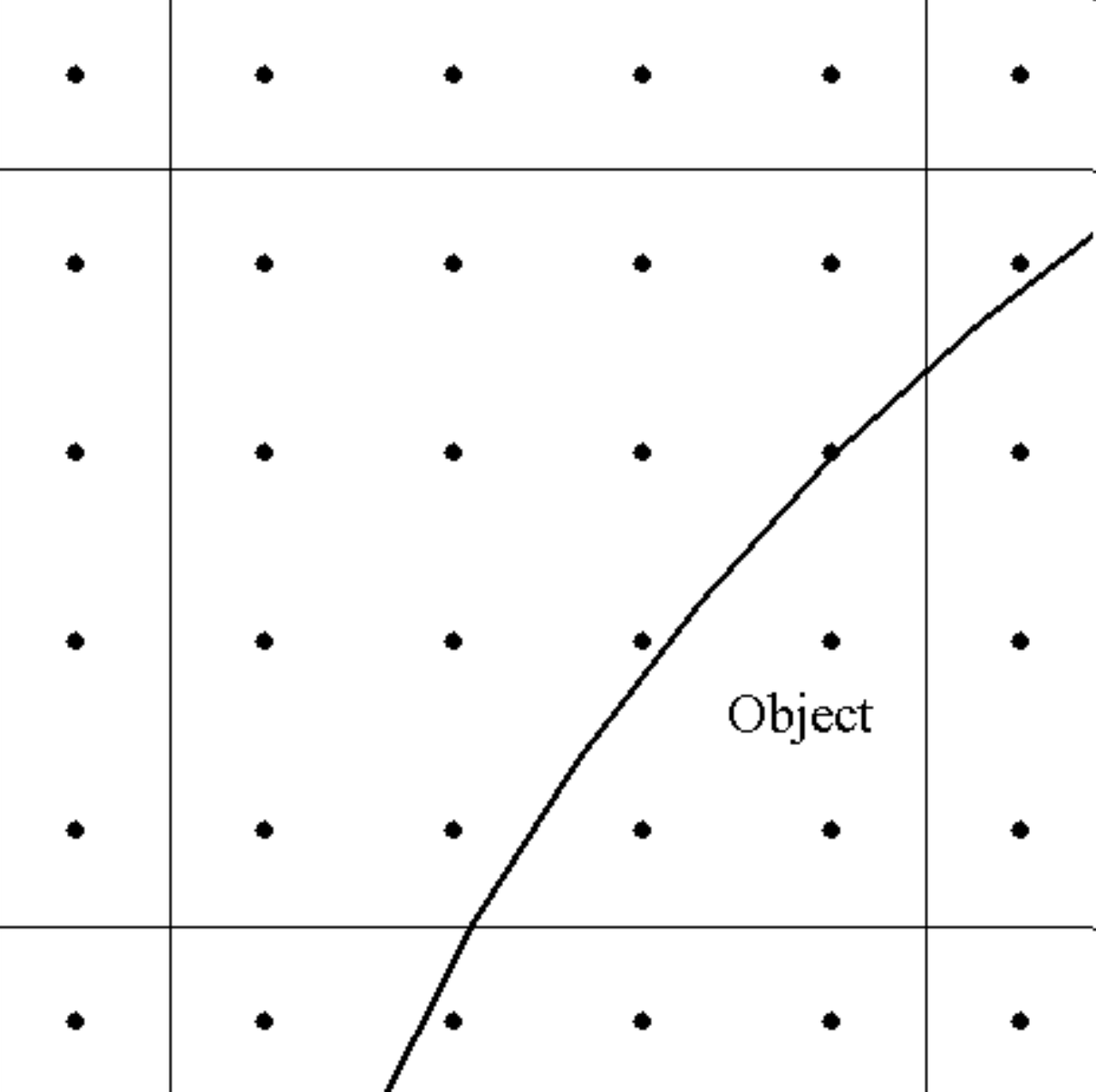} ~~~~~~
    \includegraphics[angle=0, width=1.6in]{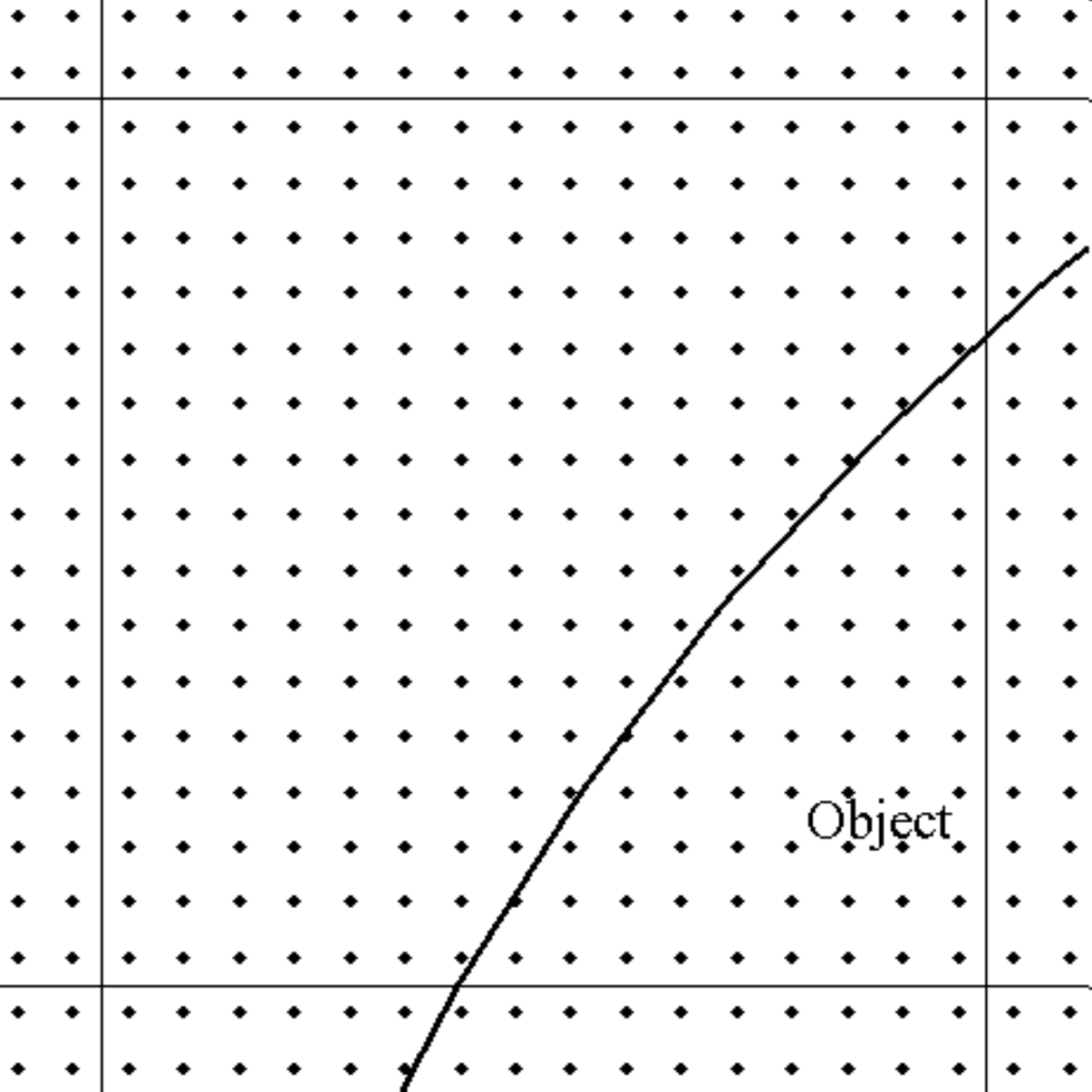}
    \caption{Effect of different particle numbers on charge distribution in interface element.}
    \label{fig: number of particles}
\end{figure}

Figure \ref{fig: number of particles} demonstrate an interface element  with $1$, or $16$, or $256$ equally distributed particles, respectively. As shown in Figure \ref{fig: number of particles}, if the number of particles is small, the charge distribution will rely more on the random fall of individual particles. On the other hand, if the number of particles is large, the effect of a single particle could be neglected. This is consistent with the actual situation.

Next, we compute the average density for different number of particles with traditional and improved PIC interpolations. We use $N$ to denote the number of particles in each interface element. We let $\mathbf{x}_i$, $i=1,\cdots N_I$ be the nodes that belongs to at least one interface element. The mean value of density at interface nodes is computed as follows
\begin{equation}\label{eq: density}
     \overline{\rho}=\frac{1}{N_{I}}\sum\limits_{i=1}^{N_{I}} \rho_h(\mathbf{x}_i),
\end{equation}
where the $\rho_h(\mathbf{x}_i)$ is the interpolated charging density at the point $\mathbf{x}_i$.
The percentage error of the density is defined by
\begin{equation}\label{eq: error density}
     E_{\rho}=\frac{|\rho- \overline{\rho}|}{|\rho|}.
\end{equation}
The $L^2$ norm error of the potential function is defined by
\begin{equation}\label{eq: error solution}
     E_{u}=\|u-u_h\| = \left(\int_\Omega |u-u_h|^2d\mathbf{x}\right)^{\frac{1}{2}}.
\end{equation}

\begin{table}[!hbt]
\newcommand{\zz}[2]{\rule[#1]{0pt}{#2}}
\newcommand{\tabincell}[2]{\begin{tabular}{@{}#1@{}}#2\end{tabular}}
\begin{center}
\begin{tabular}{|c|c|c|c|c|}
\hline
\multirow{2}{*}{$N$} & \multicolumn{2}{|c|}{Standard PIC} & \multicolumn{2}{|c|}{Improved PIC} \\
\cline{2-5}
& $\overline{\rho}$ & $E_{{\rho}}$ & $\overline{\rho}$ & $E_{{\rho}}$ \\
\hline \zz{-6pt}{12pt} $1$ & $-2.600000$ & $35.00\%$ & $-4.083332$ & $2.08\%$ \\
\hline \zz{-6pt}{12pt} $4$ & $-3.250000$ & $18.75\%$ & $-4.145832$ & $3.65\%$ \\
\hline \zz{-6pt}{12pt} $16$ & $-3.240625$ & $18.98\%$ & $-4.286457$ & $7.16\%$ \\
\hline \zz{-6pt}{12pt} $64$ & $-3.481812$ & $12.95\%$ & $-4.238281$ & $5.96\%$ \\
\hline \zz{-6pt}{12pt} $256$ & $-3.479996$ & $13.00\%$ & $-4.236654$ & $5.92\%$ \\
\hline \zz{-6pt}{12pt} $1024$ & $-3.482821$ & $12.93\%$ & $-4.240478$ & $6.01\%$ \\
\hline
\end{tabular}
\end{center}
\caption{Comparison of the error of charge distribution.}
\label{table: density}
\end{table}

\begin{table}[!hbt]
\newcommand{\zz}[2]{\rule[#1]{0pt}{#2}}
\newcommand{\tabincell}[2]{\begin{tabular}{@{}#1@{}}#2\end{tabular}}
\begin{center}
\begin{tabular}{|c|c|c|}
\hline
\zz{-6pt}{15pt} {$N$} & {IFE-PIC} & {Improved IFE-PIC} \\
\hline \zz{-6pt}{12pt}  $1$  &  $1.868830\times10^{-3}$  &  $5.160630\times10^{-4}$  \\
\hline \zz{-6pt}{12pt}  $4$  &  $1.190797\times10^{-3}$  &  $4.851096\times10^{-4}$  \\
\hline \zz{-6pt}{12pt}  $16$  &  $1.208151\times10^{-3}$  &  $6.142589\times10^{-4}$  \\
\hline \zz{-6pt}{12pt}  $64$  &  $1.281327\times10^{-3}$  &  $5.380597\times10^{-4}$  \\
\hline \zz{-6pt}{12pt}  $256$  &  $1.284765\times10^{-3}$  &  $5.359054\times10^{-4}$  \\
\hline \zz{-6pt}{12pt}  $1024$  &  $1.280203\times10^{-3}$  &  $5.403805\times10^{-4}$  \\
\hline
\end{tabular}
\end{center}
\caption{Comparison of the $L^2$ norm error of electric potential.}
\label{table: L2 solution}
\end{table}

Table \ref{table: density} and Table \ref{table: L2 solution} report the error of charging density and error of electric potential using two IFE-PIC methods. These tests are conducted on a uniform $40\times 40$ Cartesian mesh. Table \ref{table: density} clearly indicates that as we increase the number of particles belong each element, the traditional PIC interpolation become more accurate, this indicates that the traditional approach is very sensitive to the number of particles. On the other hand, the new interpolation technique is much more accurate than the traditional approach, and are more robust with respect to the number of particles. Table \ref{table: L2 solution} shows that the improved IFE-PIC scheme is about three times more accurate than the classical scheme in the $L^2$ norm of the electric potential.
\begin{remark}
Based on our current numerical experiments, simply improving only one part of IFE solver or PIC interpolation is less significant, and the behavior on accuracy varies in different choices of interface geometry, conductivity coefficient, and the number of particles.
\end{remark}

\subsection{Convergence of Potential on a Sequence of Meshes}
In this example, we test convergence of potential on a family of uniform meshes with the same number of particles in each element. We consider the number of particles is the same as in Section 4.1. Table \ref{table: diff mesh} reports the $L^2$ norm error of the potential. From the comparison, we can see that our improved IFE-PIC scheme is much more accurate on every mesh level than the widely used scheme.

\begin{table}[!hbt]
\newcommand{\zz}[2]{\rule[#1]{0pt}{#2}}
\newcommand{\tabincell}[2]{\begin{tabular}{@{}#1@{}}#2\end{tabular}}
\begin{center}
\begin{tabular}{|c|c|c|}
\hline
\zz{-6pt}{20pt} {Mesh} & {Traditional IFE-PIC} & {Improved IFE-PIC} \\
\hline \zz{-6pt}{15pt}  $10\times10$  &  $1.112170\times10^{-2}$  &  $8.131783\times10^{-3}$  \\
\hline \zz{-6pt}{15pt}  $20\times20$  &  $3.953483\times10^{-3}$  &  $2.520218\times10^{-3}$  \\
\hline \zz{-6pt}{15pt}  $40\times40$  &  $1.280203\times10^{-3}$  &  $5.403805\times10^{-4}$  \\
\hline \zz{-6pt}{15pt}  $80\times80$  &  $5.275457\times10^{-4}$  &  $1.784664\times10^{-4}$  \\
\hline \zz{-6pt}{15pt}  $160\times160$  &  $3.553394\times10^{-4}$  &  $1.313621\times10^{-4}$  \\
\hline \zz{-6pt}{15pt}  $320\times320$  &  $1.372616\times10^{-4}$  &  $1.103237\times10^{-5}$  \\
\hline \zz{-6pt}{15pt}  rate  &  $1.240231$  &  $1.771787$  \\
\hline
\end{tabular}
\end{center}
\caption{The $L^2$ norm errors of electric potential for different mesh.}
\label{table: diff mesh}
\end{table}

\section{Conclusions}
In this paper, we proposed a new IFE-PIC method for plasma simulation. The new method has improvement in both IFE solver and PIC interpolation. One prominent feature is that our new PIC interpolation has the charge-conservation property. Moreover, the improved IFE solver with partial penalty terms produce more accurate approximation around interface.

\section*{References}

\bibliography{mybibfile}

\begin{thebibliography}{10}

\bibitem{2015AdjeridChaabaneLin}
S.~Adjerid, N.~Chaabane, and T.~Lin.
\newblock An immersed discontinuous finite element method for stokes interface
  problems.
\newblock {\em Comput. Methods Appl. Mech. Engrg.}, 293:170--190, 2015.

\bibitem{PCBirch}
P.~Birch and S.~Chapman.
\newblock Two dimensional particle-in-cell simulations of the lunar wake.
\newblock {\em Physics of Plasma}, 9(5):1785--1789, 2002.

\bibitem{2015CaoChuWangCaoXiaZhang}
H.~Cao, Y.~Chu, E.~Wang, Y.~Cao, G.~Xia, and Z.~Zhang.
\newblock Numerical simulation study on barrel erosion of ion thruster
  accelerator grid.
\newblock {\em Journal of Propulsion and Power}, 31(6):1785--1792, 2015.

\bibitem{2015CaoLiShanCaoZheng}
H.~Cao, Q.~Li, K.~Shan, Y.~Cao, and L.~Zheng.
\newblock Effect of preionization on the erosion of the discharge channel wall
  in a hall thruster using a kinetic simulation.
\newblock {\em IEEE Transactions on Plasma Science}, 43(1):130--140, 2015.

\bibitem{YCao_YChu_XMHe_TLin}
Y.~Cao, Y.~Chu, X.~He, and T.~Lin.
\newblock An iterative immersed finite element method for an electric potential
  interface problem based on given surface electric quantity.
\newblock {\em Journal of Computational Physics}, 281(15):82--95, 2015.

\bibitem{2016CaoChuZhangZhang}
Y.~Cao, Y.~Chu, X.~Zhang, and X.~Zhang.
\newblock Immersed finite element methods for unbounded interface problems with
  periodic structures.
\newblock {\em J. Comput. Appl. Math.}, 307:72--81, 2016.

\bibitem{2011ChuCaoHeLuo}
Y.~Chu, Y.~Cao, X.~He, and M.~Luo.
\newblock Asymptotic boundary conditions with immersed finite elements for
  interface magnetostatic/electrostatic field problems with open boundary.
\newblock {\em Comput. Phys. Comm.}, 182(11):2331--2338, 2011.

\bibitem{2017ChuHanCaoHeWang}
Y.~Chu, D.~Han, Y.~Cao, X.~He, and J.~Wang.
\newblock An immersed-finite-element particle-in-cell simulation tool for
  plasma surface interaction.
\newblock {\em International Journal of Numerical Analysis and Modeling},
  14(2):175--200, 2017.

\bibitem{JDawson}
J.~Dawson.
\newblock One-dimensional plasma model.
\newblock {\em The Physics of Fluids}, 5(4):445--459, 1962.

\bibitem{OCEldridge}
O.~Eldridge and M.~Feix.
\newblock One-dimensional plasma model at thermodynamic equilibrium.
\newblock {\em The Physics of Fluids}, 5(9):1076--1080, 1962.

\bibitem{2010GongLi}
Y.~Gong and Z.~Li.
\newblock Immersed interface finite element methods for elasticity interface
  problems with non-homogeneous jump conditions.
\newblock {\em Numer. Math. Theory Methods Appl.}, 3(1):23--39, 2010.

\bibitem{2016HanWangHe}
D.~Han, J.~Wang, and X.~He.
\newblock a nonhomogeneous immersed finite element particle in cell method for
  modeling dielectric surface charging in plasmas.
\newblock {\em IEEE Transactions on Plasma Science}, 44(8):1326--1332, 2016.

\bibitem{DHan_PWang}
D.~Han, P.~Wang, X.~He, T.~Lin, and J.~Wang.
\newblock A 3{D} immersed finite element method with non-homogeneous interface
  flux jump for applications in particle-in-cell simulations of plasma-lunar
  surface interactions.
\newblock {\em Journal of Computational Physics}, 321:965--980, 2016.

\bibitem{2008HeLinLin}
X.~He, T.~Lin, and Y.~Lin.
\newblock Approximation capability of a bilinear immersed finite element space.
\newblock {\em Numer. Methods Partial Differential Equations},
  24(5):1265--1300, 2008.

\bibitem{2011HeLinLin}
X.~He, T.~Lin, and Y.~Lin.
\newblock Immersed finite element methods for elliptic interface problems with
  non-homogeneous jump conditions.
\newblock {\em Int. J. Numer. Anal. Model.}, 8(2):284--301, 2011.

\bibitem{2013HeLinLinZhang}
X.~He, T.~Lin, Y.~Lin, and X.~Zhang.
\newblock Immersed finite element methods for parabolic equations with moving
  interface.
\newblock {\em Numer. Methods Partial Differential Equations}, 29(2):619--646,
  2013.

\bibitem{1981HockneyEastwood}
R.~W. Hockney and J.~W. Eastwood.
\newblock {\em Computer simulation using particles}.
\newblock Taylor and Francis, 1981.

\bibitem{2015JianChuCaoCaoHeXia}
H.~Jian, Y.~Chu, H.~Cao, Y.~Cao, X.~He, and G.~Xia.
\newblock Three-dimensional {IFE-PIC} numerical simulation of background
  pressure's effect on accelerator grid impingement current for ion optics.
\newblock {\em Vacuum}, 116:130--138, 2015.

\bibitem{2005KafafyLinLinWang}
R.~Kafafy, T.~Lin, Y.~Lin, and J.~Wang.
\newblock Three-dimensional immersed finite element methods for electric field
  simulation in composite materials.
\newblock {\em Internat. J. Numer. Methods Engrg.}, 64(7):940--972, 2005.

\bibitem{2005KafafyWangLin}
R.~Kafafy, J.~Wang, and T.~Lin.
\newblock A hybrid-grid immersed-finite-element particle-in-cell simulation
  model of ion optics plasma dynamics.
\newblock {\em Dyn. Contin. Discrete Impuls. Syst. Ser. B Appl. Algorithms},
  12(Suppl. Vol. 12b):1--16, 2005.

\bibitem{SKimura}
S.~Kimura and T.~Nakagawa.
\newblock Electromagnetic full particle simulation of the electric field
  structure around the moon and the lunar wake.
\newblock {\em Earth, Planets and Space}, 60(6):594--599, 2008.

\bibitem{LangdonAB_Birdsall}
A.~Langdon and C.~Birdsall.
\newblock Theory of plasma simulation using finite-size particles.
\newblock {\em The Physics of Fluids}, 13(8):2115--2122, 1970.

\bibitem{1998Li}
Z.~Li.
\newblock The immersed interface method using a finite element formulation.
\newblock {\em Appl. Numer. Math.}, 27(3):253--267, 1998.

\bibitem{2004LiLinLinRogers}
Z.~Li, T.~Lin, Y.~Lin, and R.~C. Rogers.
\newblock An immersed finite element space and its approximation capability.
\newblock {\em Numer. Methods Partial Differential Equations}, 20(3):338--367,
  2004.

\bibitem{2003LiLinWu}
Z.~Li, T.~Lin, and X.~Wu.
\newblock New {C}artesian grid methods for interface problems using the finite
  element formulation.
\newblock {\em Numer. Math.}, 96(1):61--98, 2003.

\bibitem{2013LinLinZhang2}
T.~Lin, Y.~Lin, and X.~Zhang.
\newblock Immersed finite element method of lines for moving interface problems
  with nonhomogeneous flux jump.
\newblock In {\em Recent advances in scientific computing and applications},
  volume 586 of {\em Contemp. Math.}, pages 257--265. Amer. Math. Soc.,
  Providence, RI, 2013.

\bibitem{2013LinLinZhang1}
T.~Lin, Y.~Lin, and X.~Zhang.
\newblock A method of lines based on immersed finite elements for parabolic
  moving interface problems.
\newblock {\em Adv. Appl. Math. Mech.}, 5(4):548--568, 2013.

\bibitem{2015LinLinZhang}
T.~Lin, Y.~Lin, and X.~Zhang.
\newblock Partially penalized immersed finite element methods for elliptic
  interface problems.
\newblock {\em SIAM J. Numer. Anal.}, 53(2):1121--1144, 2015.

\bibitem{2013LinSheenZhang}
T.~Lin, D.~Sheen, and X.~Zhang.
\newblock A locking-free immersed finite element method for planar elasticity
  interface problems.
\newblock {\em J. Comput. Phys.}, 247:228--247, 2013.

\bibitem{2015LinYangZhang1}
T.~Lin, Q.~Yang, and X.~Zhang.
\newblock {\it A Priori} error estimates for some discontinuous {Galerkin}
  immersed finite element methods.
\newblock {\em J. Sci. Comput.}, 65(3):875--894, 2015.

\bibitem{2012LinZhang}
T.~Lin and X.~Zhang.
\newblock Linear and bilinear immersed finite elements for planar elasticity
  interface problems.
\newblock {\em J. Comput. Appl. Math.}, 236(18):4681--4699, 2012.

\bibitem{LASchwager}
L.~Schwager and C.~Birdsall.
\newblock Collector and source sheaths of a finite ion temperature plasma.
\newblock {\em The Physics of Fluids B: Plasma Physics}, 2(5):1057--1068, 1990.

\bibitem{VlasovAA}
A.~Vlasov.
\newblock On vibration properties of electron gas.
\newblock {\em Journal of Experimental and Theoretical Physics}, 8(3):291,
  1938.

\end{thebibliography}

\end{document}